\documentclass[11pt]{amsart}
\usepackage{hyperref}
\usepackage{graphicx}
\usepackage{enumerate}
\usepackage{amsthm}
\usepackage{epsfig}

\usepackage{pgf}
\usepackage{tikz}
\usetikzlibrary{arrows,shapes,er,automata,backgrounds,topaths,trees}

\usepackage{slashbox}

\newcommand*{\CopyCounter}[2]{%
  \expandafter\def\csname c@#2\endcsname{\csname c@#1\endcsname}%
  \expandafter\def\csname p@#2\endcsname{\csname p@#1\endcsname}%
  \expandafter\def\csname the#2\endcsname{\csname the#1\endcsname}}

\numberwithin{Theorem}{section}
\CopyCounter{Theorem}{Proposition}
\CopyCounter{Theorem}{ProposedProblem}
\CopyCounter{Theorem}{Property}
\CopyCounter{Theorem}{Claim}
\CopyCounter{Theorem}{Lemma}
\CopyCounter{Theorem}{Corollary}
\CopyCounter{Theorem}{Conjecture}
\CopyCounter{Theorem}{Definition}
\CopyCounter{Theorem}{Example}
\CopyCounter{Theorem}{Remark}
\CopyCounter{Theorem}{Question}
\CopyCounter{Theorem}{Condition}
\CopyCounter{Theorem}{Criterion}
\CopyCounter{Theorem}{Observation}
\theoremstyle{plain}%% needs amsthm.sty

\newtheorem{property}[Property]{Property}

\theoremstyle{definition}%% needs amsthm.sty

\newtheorem{remark}[Remark]{Remark}

\newcommand{\bq}{\begin{equation}}
\newcommand{\eq}{\end{equation}}

\newcommand{\grad}{\nabla}

\newcommand{\bm}[1]{\mbox{\boldmath${#1}$}}

\newcommand{\domain}{\Omega}
\newcommand{\cdomain}{\bar{\Omega}}
\newcommand{\boundary}{\partial \domain}

\newcommand{\inBoundary}{\mathcal{B}_{\text{in}}}
\newcommand{\outBoundary}{\mathcal{B}_{\text{out}}}
\newcommand{\fBoundary}{\mathcal{B}_{\text{f}}}

\newcommand{\btilde}{\tilde{\bm{b}}}
\newcommand{\xtilde}{\tilde{\bm{x}}}

\newcommand{\xhat}{\hat{\bm{x}}}
\newcommand{\yhat}{\hat{\bm{y}}}

\newcommand{\x}{\bm{x}}

\newcommand{\ba}{\bm{a}}
\newcommand{\bb}{\bm{b}}
\newcommand{\bc}{\bm{c}}

\newcommand{\y}{\bm{y}}

\newcommand{\J}{{\mathcal{J}}} 
\newcommand{\Jhat}{{\mathcal{\hat{J}}}}

\newcommand{\T}{{\mathcal{T}}}

\newcommand{\R}{\bm{R}}

\newcommand{\fB}{\bm{f}}

\newcommand{\htilde}{\tilde{h}}
\newcommand{\blds}[1]{\mbox{\scriptsize \boldmath $#1$}}

\DeclareMathOperator*{\distance}{distance}

  \renewenvironment{thebibliography}[1]{%
    \begin{oldthebibliography}{#1}%
      \setlength{\parskip}{.3ex}%
      \setlength{\itemsep}{.3ex}%
  }%
  {%
    \end{oldthebibliography}%
  }

\begin{document}

\title[Efficient method for multiobjective optimal control]
{An efficient method for multiobjective optimal control 
and optimal control subject to integral constraints.}
\author{Ajeet Kumar}
\address{Department of Theoretical \& Applied Mechanics, Cornell University}
\email{ak428@cornell.edu}
\author{Alexander Vladimirsky}
\address{Department of Mathematics, Cornell University}
\email{vlad@math.cornell.edu}
\thanks{This research has been supported in part by the NSF grant DMS-0514487.}
\date{\today}

\begin{abstract}
\noindent
We introduce a new and efficient numerical method for 
multicriterion optimal control and single criterion optimal control
under integral constraints.  The approach is based on extending
the state space to include information on a ``budget'' remaining to satisfy each 
constraint;  the augmented Hamilton-Jacobi-Bellman PDE is then solved
numerically.  The efficiency of our approach hinges on the causality 
in that PDE, i.e., the monotonicity of characteristic curves in
one of the newly added dimensions.  A semi-Lagrangian ``marching'' method
is used to approximate the discontinuous viscosity solution efficiently.
We compare this to a recently introduced ``weighted sum'' based algorithm
for the same problem \cite{MitchellSastry}.  
We illustrate our method using examples from 
flight path planning and robotic navigation in the presence of friendly 
and adversarial observers.
\end{abstract}

\maketitle

\section{Introduction.}
\label{s:intro}

In the continuous setting, deterministic optimal control problems 
are often studied from the point of view of dynamic programming; see, 
e.g., \cite{BardiDolcetta}, \cite{BressanPiccoli}.
A choice of the particular control $\ba(t)$ determines the trajectory
$\y(t)$ in the space of system-states $\domain \subset \R^n$.
A running cost $K$ is integrated along that trajectory and the terminal cost $q$
is added, yielding the total cost associated with this control.
A {\em value function} $u$, describing the minimum cost to pay starting
from each system state, is shown to be the unique viscosity solution of
the corresponding Hamilton-Jacobi PDE.  Once the value function has been
computed, it can be used to approximate optimal feedback control.
We provide an overview of this classic approach in section \ref{s:traditional}.

However, in realistic applications practitioners usually need to optimize
by many different criteria simultaneously.  
For example, given a vehicle starting at $\x \in \domain$ and trying 
to ``optimally'' reach some target $\T$, the above framework allows to find 
the most fuel efficient trajectories and the fastest trajectories, but these
generally will not be the same.  
A natural first step is to compute the total time taken along 
the most fuel-efficient trajectory and the total amount of fuel needed to follow
the fastest trajectory.  Computational efficiency requires a method for computing 
this simultaneously for all starting states $\x$.  A PDE-based approach for this task
is described in section \ref{ss:secondary-along-primary}.  

This, however, does not yield answers to two more practical questions: 
what is the fastest trajectory from $\x$ to $\T$ without using more 
than the specified amount of fuel?  Alternatively, 
what is the most fuel-efficient trajectory from $\x$, provided 
the vehicle has to reach $\T$ by the specified time?
We will refer to such trajectories as {\em constrained-optimal}.

One approach to this more difficult problem is the {\em Pareto optimization}:
finding a set of trajectories, which are optimal in the sense that no improvement in 
fuel-efficiency is possible without spending more time (or vice versa).
This defines a {\em Pareto front} -- a curve in a time-fuel plane, 
where each point corresponds to time \& fuel needed along some Pareto-optimal 
trajectory.  This approach is generally computationally expensive, especially
if a Pareto front has to be found for each starting state $\x$ separately.
The current state of the art for this problem has been developed by Mitchell and Sastry
\cite{MitchellSastry} and described in section \ref{ss:Mitchell-Sastry}.
Their method is based on the usual {\it weighted sum} approach 
to multiobjective optimization \cite{Marler}.  
A new running cost $K$ is defined as a weighted average of several competing 
running costs $K_i$'s, and the corresponding Hamilton-Jacobi PDE 
is then solved to obtain one point on the Pareto front.  
The coefficients in the weighted sum are then varied and the process is 
repeated until a solution satisfying all constraints is finally found.  
Aside from the computational cost, 
the obvious disadvantage of this approach is that only a convex part of 
the Pareto front can be obtained by weighted sum methods \cite{Das}, which
may result in selecting suboptimal trajectories.  In addition,
recovering the entire Pareto front for each $\x \in \domain$
is excessive and unnecessary when the real goal is to solve
the problem for a fixed list of constraints (e.g.,
maximum fuel or maximum time available).

Our own approach (described in section \ref{ss:new-approach}) 
remedies these problems by systematically constructing
the exact portion of Pareto front relevant to the above constraints  
{\em for all $\x \in \domain$ simultaneously}.
Given $\domain \subset \R^n$ and $r$ additional integral constraints,
we accomplish this by solving a single augmented partial differential equation on
a $(r+n)$-dimensional domain.  Our method has two key advantages.
First, it does not rely on any assumptions about the convexity of Pareto front.
Secondly, the PDE we derive has a special structure, allowing for a very efficient 
marching method.  Our approach can be viewed as a generalization of the classic 
equivalency of Bolza and Mayer problems \cite{BressanPiccoli}.
The idea of accommodating integral constraints by extending the state space
is not new.  It was previously used by Isaacs to derive the properties of
constrained-optimal strategies for differential games \cite{IsaacsBook}.  
More recently,
it was also used in infinite-horizon control problems 
by Soravia \cite{Soravia_IntegralConstraint} and 
Motta \& Rampazzo \cite{MottaRampazzo}
to prove the uniqueness of the (lower semi-continuous) 
viscosity solution to the augmented PDE.
However, the above works explored the theoretical issues only and,
to the best of our knowledge, ours is the first practical numerical method based on
this approach.  In addition, we also show the relationship between
optimality under constraints and Pareto optimality for feasible trajectories.

The computational efficiency of our method is deeply related to the general difference in 
numerical methods for time-dependent and static first-order equations.  
In optimal control problems, time-dependent HJB PDEs result from 
{\em finite-horizon} problems or problems with time-dependent dynamics and running cost.
Static HJB PDEs usually result from {\em exit-time} or {\em infinite-horizon} 
problems with time-independent (though perhaps time-discounted) dynamics 
and running cost.  In the time-dependent case, efficient numerical methods are
typically based on time-marching.  In the static case, a naive approach
involves iterative solving of a system of discretized equations.  Several popular
approaches were developed precisely to avoid these iterations either by 
space-marching (e.g., \cite{Tsitsiklis_conference, SethFastMarcLeveSet, SethVlad2}), 
or by embedding into a higher-dimensional time-dependent problem 
(via Level Set Methods, e.g., \cite{Osher_Dirichlet}), or by 
treating one of the spatial directions as if it were time (resulting
in a ``paraxial'' approximation; see, e.g., \cite{QianSymes1}).  
For reader's convenience,
we provide a brief overview of these approaches in sections \ref{ss:causal_discr}
and \ref{ss:marching}.
Our key observation is that the augmented ``static'' PDE has {\em explicit causality},
allowing simple marching (similar to time-marching) in the secondary cost direction.

Our semi-Lagrangian method is described in section \ref{ss:new-numerics}.
Since the augmented PDE is solved on a higher-dimensional domain, any restriction of 
that domain leads to substantial computational savings.
In section \ref{ss:reduced_domain} we explain how this can be accomplished by solving
static PDEs in $\domain$ for each individual cost.  This additional step also 
yields improved boundary conditions for the primary value function in $\R^{n+r}$.

In section \ref{s:numerical-results} we illustrate our method using examples 
from robotic navigation (finding shortest/quickest paths, while avoiding (or seeking) exposure to
stationary observers) and a test-problem introduced in
\cite{MitchellSastry} (planning a flight-path for an airplane 
to minimize the risk of encountering a storm while satisfying constraints 
on fuel consumption).  
Finally, in section \ref{s:conclusions} we discuss the limitations
of our approach and list several directions for future research.

\section{Single-criterion dynamic programming.}
\label{s:traditional}

\subsection{Exit-time optimal control.}
\label{ss:static}

To begin, we consider a general deterministic exit-time optimal control problem.
This is a classic problem and our discussion follows the description in \cite{BardiDolcetta}.
Suppose $\domain \subset \R^n$ is an open bounded set of all possible ``non-terminal''
states of the system, while $\boundary$ is the set of all terminal states.
For every starting state $\x \in \domain$,
the goal is to find the cheapest way to leave $\domain$. 

Suppose $A \in \R^m$ is a compact set of control values,
and the set of admissible controls $\mathcal{A}$ consists of all measurable functions 
$\ba: \R \mapsto A$.  The dynamics of the system is defined by 
\begin{eqnarray}
\nonumber
\y^{\prime}(t) &=& \fB(\y(t), \ba(t)),\\
\y (0) &=& \x \in \domain,
\label{eq:auton_dynamics}
\end{eqnarray}
where $\y(t)$ is the system state at the time $t$,
$\x$ is the initial system state, 
and $\fB: \cdomain \times A \mapsto \R^n$ is the velocity.
The exit time associated with this control is 
\begin{equation}
T_{\x, \ba} = \min \{ t \in \R_{+,0} | \y(t) \in \boundary \}.
\label{eq:exit-time}
\end{equation}
The problem description also includes 
the running cost $K: \cdomain \times A \mapsto \R$ and
the terminal (non-negative, possibly infinite) cost 
$q: \boundary \mapsto \left(\R_{+,0} \cup \{+\infty\} \right)$.  
This allows to specify the total cost of using the control $\ba(\cdot)$ starting from $\x$:
$$
\J(\x, \ba(\cdot)) = \int_0^{T_{\x,\ba}} K (\y(t), \ba(t)) \, dt 
\, + \, q(\y(T_{\x,\ba})).
$$

We will make the following assumptions throughout the rest of this paper:
\begin{enumerate}[{\quad (A}1)]
\setcounter{enumi}{-1}
\item 
Function $q$ is lower semi-continuous and $\min_{\boundary} q < +\infty$.
\item 
Functions $\fB$ and $K$ are Lipschitz-continuous.
\item
There exist constants $k_1, k_2$ such that
$0 < k_1 \leq K(\x, \ba) \leq k_2$ for $\forall \x \in \cdomain, \ba \in A.$
\item
For every $\x \in \domain$, the {\em scaled vectogram} 
$$V(\x) = \{ \fB(\x, \ba) / K(\x, \ba) \, \mid \, \ba \in A \}$$
is a strictly convex set, containing the origin in its interior.
\end{enumerate}

The key idea of {\em dynamic programming} is to introduce the value function
$u(\x)$, describing the minimum cost needed to exit $\domain$ starting from $\x$:
\begin{equation}
u(\x) \; = \; \inf_{\ba(\cdot) \in \mathcal{A}} \J(\x, \ba(\cdot)).
\label{eq:value}
\end{equation}
Bellman's optimality principle \cite{Bellman_dynamic} shows that,
for every sufficiently small $\tau > 0$,
\begin{equation}
u(\x) \; = \; \inf_{\ba(\cdot)} \left\{ \int_0^{\tau} K (\y(t), \ba(t)) \, dt 
\; + \; u(\y(\tau)) \right\},
\label{eq::Bellman_contin}
\end{equation}
where $\y(\cdot)$ is a trajectory corresponding to a chosen control $\ba(\cdot)$.
If $u(\x)$ is smooth, a Taylor expansion of the above formula can be used to formally 
show that $u$ is the solution of a static Hamilton-Jacobi-Bellman PDE:
\begin{eqnarray}
\nonumber
\min\limits_{\ba \in A}
\left\{ 
K(\x, \ba) + \nabla u (\x) \cdot \fB ( \x, \ba )
\right\}
\, = \, 0,
&& \text{ for $\x \in \domain$};\\
u(\x) \, = \, q(\x), && \text{ for $\x \in \boundary$}.
\label{eq:HJB_general}
\end{eqnarray}
Unfortunately, a smooth solution to Eqn. \ref{eq:HJB_general} 
might not exist even for smooth $\fB$,$K$, $q$, and $\boundary$.  
Generally, this equation has infinitely 
many weak Lipschitz-continuous solutions, but the unique {\em viscosity solution}
can be defined using additional conditions on smooth test functions \cite{CranLion, CranEvanLion}.
It is a classic result that the viscosity solution of this PDE 
coincides with the value function of the above control problem.

Under the above assumptions 
the value function $u(\x)$ is locally Lipschitz-continuous on $\domain$, 
an optimal control $\ba(\cdot)$ actually exists for every $\x \in \domain$
(i.e., $\bm{\min}$ can be used instead of $\bm{\inf}$ in formula \ref{eq:value}),
and the minimizing control value $\ba$ (in equation \ref{eq:HJB_general})
is unique wherever $\nabla u(\x)$ is defined \cite{BardiDolcetta}.
The characteristic curves of this PDE are the optimal trajectories 
of the control problem.  The points where $\nabla u$ is undefined 
are precisely those, where multiple characteristics intersect
(or, alternatively, 
the points for which multiple optimal trajectories are available).
We will let $\Gamma$ be a set of all such points.
By Rademacher's theorem, $\Gamma$ has a measure zero in $\cdomain$.

We note that the assumptions (A2-A3) can be relaxed at the cost of
additional technical details.
For example, if $V(\x)$ is not convex, the existence of an optimal 
control is not guaranteed even though the value function can still
be recovered from the PDE~\eqref{eq:HJB_general} and there exist
suboptimal controls whose cost is arbitrarily close to $u(\x)$.
On the other hand, if $V(\x)$ does not contain the origin in its interior,
then the {\em reachable set} 
$$
\mathcal{R} = \{ \x \in \domain \, \mid \, 
\text{ there exists a control leading from $\x$ to $\boundary$ in finite time}
\}
$$
need not be the entire $\domain$.  In that case, 
$\partial \mathcal{R}$ is a {\em free boundary}.
We refer readers to \cite{BardiDolcetta} for further details.

The above framework is also flexible enough to
describe the task of optimally reaching some compact target 
$\T \subset \domain$ without leaving $\domain$.
To do that, we can simply pose the problem on a new 
domain $\domain^{new} = \domain \backslash \T$, defining
the new exit cost $q^{new}$ to be infinite on $\boundary$
and finite on $\partial \T$.
The assumptions (A0-A3) guarantee 
the continuity of the value function on $\domain$ even in the presence
of state-constraints; $k_1 > 0$ and the fact that the origin is in the interior
of $V(\x)$ yield both
Soner's ``inward pointing'' condition along the boundary of the constraint set 
(as in \cite{Soner})
and the local controllability 
(as in \cite{BardiFalcone}, for example).

\begin{remark}
The continuity of $u$ on $\cdomain$ is a more subtle issue 
requiring additional {\em compatibility conditions} on $q$ even if 
that function is continuous;
otherwise, the boundary conditions are satisfied by the value function
``in viscosity sense'' only \cite{BardiDolcetta}.
However, due to our very strong controllability assumption (A3), 
the local Lipschitz-continuity 
of $u$ in the interior is easy to show even if 
$q$ is discontinuous, as allowed by (A0).  
Without assumption (A3) or its equivalent,
discontinuous boundary data typically leads to discontinuities 
in the value function in the interior as well.  
Such is the case for the augmented PDE defined in section \ref{ss:new-approach}.
\label{rem:boundary_continuity}
\end{remark}

Before continuing with the general case we consider two particularly
useful subclasses of problems.

{\em Geometric dynamics.}
Suppose $A =\{ \ba \in \R^n \mid |\ba| \leq 1\}$ and, for all 
$\x \in \domain, \ba \in A \backslash \{0\}$,
$$
K(\x, \ba) \geq  |\ba| K(\x, \ba/|\ba|), \quad
\fB ( \x, \ba ) = f( \x, \ba/|\ba|) \ba, \quad
\text{ and }
\fB ( \x, 0 ) = 0. 
$$ 
Then it is easy to show that the cost of any trajectory
is reduced by traversing it as quickly as possible; i.e.,
we can redefine $A = S_{n-1}=\{ \ba \in \R^n \mid |\ba| = 1\}$ without 
affecting the value function.
In that case, the control value $\ba$ is simply our choice for the direction 
of motion and $f$ is the speed of motion in that direction.
The equation~\eqref{eq:HJB_general} now becomes
\begin{equation}
\min\limits_{\ba \in S_{n-1}}
\left\{ 
K(\x, \ba) + \left(\nabla u (\x) \cdot \ba \right) f( \x, \ba )
\right\}
\, = \, 0.
\label{eq:HJB_geometric}
\end{equation}
A further simplification is possible if the speed and cost are isotropic, i.e.,
$f( \x, \ba ) = f(\x)$ and $K( \x, \ba ) = K(\x)$.  
In this case, the minimizing control value is
$\ba = -\nabla u (\x) / |\nabla u (\x)|$ and 
~\eqref{eq:HJB_general} reduces to an {\em Eikonal PDE}:
\begin{equation}
|\nabla u (\x)| \, f( \x ) \; = \; K(\x).
\label{eq:Eikonal}
\end{equation}
The characteristics of the Eikonal equation coincide with the gradient lines
of its viscosity solution.  
That special property can be used to construct particularly efficient numerical
methods for this PDE; see section \ref{ss:marching}.

{\em Time-optimal control.}
If $K(\x,\ba) = 1$ for all $\x$ and $\ba$, then 
$
\J(\x, \ba(\cdot)) = T_{\x,\ba} \, + \, q(\y(T_{\x,\ba})).
$
Interpreting $q$ as an {\em exit time penalty},
we can restate this as a problem of finding a {\em time-optimal control}
starting from $\x$.  The PDE~\eqref{eq:HJB_general} becomes
\begin{equation}
\max\limits_{\ba \in A}
\left\{ 
-\nabla u (\x) \cdot \fB ( \x, \ba )
\right\}
\, = \, 1.
\label{eq:HJB_time-optim}
\end{equation}
We note that the value function for {\em every} exit-time optimal control problem 
satisfying assumptions (A0-A3)
can be reduced to this time-optimal control case by setting 
$K^{new} = 1$ and $\fB^{new} ( \x, \ba ) = \fB ( \x, \ba ) / K ( \x, \ba )$;
a proof of this for the case of geometric dynamics can be found in \cite{VladThesis}.

A combination of these two special cases is attained when $K=1$ and $f=1$,
leading to a PDE $|\nabla u(\x)|=1$.  If the boundary condition is $q=0$,
then $u(\x)$ is simply the distance from $\x$ to $\boundary$.

\subsection{Semi-Lagrangian discretizations.}
\label{ss:discr}

Suppose $X$ is a grid or a mesh on the domain $\cdomain$
and for simplicity assume that $\boundary$ is well-discretized by $\partial X$.
We will use $M=|X|$ to denote the total number of meshpoints in $X$.

A natural first-order accurate semi-Lagrangian discretization of equation 
~\eqref{eq:HJB_general} is obtained by assuming that the control value
is held fixed for some small time $\tau$. If $U(\x)$ is the approximation
of the value function at the mesh point $\x \in X$, then
the optimality principle yields the following system:
\begin{eqnarray}
\nonumber
U(\x) &=& \min\limits_{\ba \in A}
\left\{
\tau K(\x, \ba) \; + \; U
\left( 
\x + \tau \fB(\x, \ba)
\right)
\right\},
\qquad
\forall \x \in X \backslash \boundary;\\
\label{eq:Falcone-style}
\\
\nonumber
U(\x) &=& q(\x), \qquad 
\forall \x \in \partial X.
\end{eqnarray}
Since $\xtilde_{\blds{a}} = \x + \tau \fB(\x, \ba)$ usually is not a gridpoint,
a first-order accurate reconstruction is needed to approximate $U(\xtilde)$
based on values of $U$ at nearby meshpoints\footnote{
If $\x$ is close to $\boundary$, it is possible that 
$\xtilde_{\blds{a}} \not \in \cdomain$.  This can be handled by either enlarging
$X$ to cover some neighborhood of $\domain$ (see, e.g., \cite{BardiFalconeSoravia1})
or by decreasing $\tau$ in such cases to make $\xtilde_{\blds{a}} \in \boundary$.  
The latter strategy is also adopted in our implementation of the method in section \ref{ss:new-numerics}.}.
Discretizations of this type were introduced by Falcone 
for time-discounted infinite-horizon problems, which can be related 
to the above exit-time problem using the Kruzhkov transform
\cite{Falcone_InfHor, Falc_Dial}.

In case of geometric dynamics, another natural discretization is based
on the assumption that the direction of motion is held constant until reaching
a boundary of a simplex.  For notational simplicity, suppose that $n=2$
and $S(\x)$ is the set of all triangles in $X$ with a vertex at $\x$.
For every $s \in S(x)$, denote its other vertices as $\x_{s1}$ and $\x_{s2}$.
Suppose $\xtilde_{\blds{a}} = \x + \tau_{\blds{a}} f(\x, \ba) \ba$ lies on the segment
$\x_{s1} \x_{s2}$.  Let 
$\Xi = \left\{ \xi = (\xi_1, \xi_2) \, \mid \, 
\xi_1 + \xi_2 = 1 \text { and } \xi_1, \xi_2 \geq 0 \right\}.$
Then $\xtilde_{\blds{a}} = \xi_1 \x_{s1} + \xi_2 \x_{s2}$ for some $\xi \in \Xi$
and $U(\xtilde_{\blds{a}}) = \xi_1 U(\x_{s1}) + \xi_2 U(\x_{s2})$.
Alternatively, given $\xtilde_{\xi} = \xi_1 \x_{s1} + \xi_2 \x_{s2}$,
we can define $\ba_{\xi} = (\xtilde_{\xi} - \x) / |\xtilde_{\xi} - \x|$
and $\tau_{\xi} = |\xtilde_{\xi} - \x| / f(\x, \ba_{\xi})$.
This yields the following system of discretized equations:
\begin{eqnarray}
\nonumber
U(\x) & = &
\min\limits_{s \in S(\x)} \,
\min\limits_{\xi \in \Xi}
 \left\{ 
\tau_{\xi} K(\x, \ba_{\xi})
\, + \, \left( \xi_1 U(\x_{s1}) + \xi_2 U(\x_{s2}) \right)
\right\},  \qquad \forall \x \in X \backslash \boundary;\\
\nonumber
U(\x) &=& q(\x),
  \qquad \forall \x \in \partial X.\\
\label{eq:geometric_aniso}
\end{eqnarray}
Discretizations of this type were used by Gonzales
and Rofman in \cite{GonzalezRofman}.  Both discretizations
~\eqref{eq:Falcone-style} and 
~\eqref{eq:geometric_aniso} were also earlier derived
by Kushner and Dupuis approximating continuous optimal control 
by controlled Markov processes \cite{KushnerDupuis}.
In the appendix of \cite{SethVlad3} we demonstrated that 
~\eqref{eq:geometric_aniso} is also equivalent to a first-order
upwind finite difference approximation on the same mesh $X$.

\subsection{Causality, dimensionality \& computational efficiency.}
\label{ss:causal_discr}

We note that both of the above discretizations result in 
a system of $M$ non-linear coupled equations.
Finding a numerical solution to that system efficiently is
an important practical problem in itself.  

Suppose all meshpoints in $X$ are ordered $\x_1, \ldots, \x_M$ and
$U = (U_1, \ldots, U_M)$ is a vector of corresponding approximate values.
The above discretized equations can be formally
written as $U=\mathcal{F}(U)$ and one simple approach is 
to recover $U$ by fixed-point iterations, i.e., 
$U^{k+1} = \mathcal{F}(U^k)$, where $U^0$ is 
an appropriate initial guess for $U$.
This procedure is generally quite inefficient
since each iteration has a $O(M)$ computational cost
and a number of iterations needed is at least $O(M^{1/n})$.

A Gauss-Seidel relaxation of the above is a typical
practical modification, where the entries of $U^{k+1}$
are computed sequentially and 
the new values are used as soon as they become available:
$U_i^{k+1} = \mathcal{F}_i (U_1^{k+1}, \ldots, U_{i-1}^{k+1}, U_i^k, \ldots, U_M^k).$
The number of iterations required to converge will now heavily
depend on the PDE, the particular discretization and the ordering of the meshpoints.
We will say that a discretization is {\em causal}
if there exists an ordering of meshpoints such that
the convergence is attained after a {\em single} Gauss-Seidel iteration.
For example, if the dynamics of the control problem is such that
one of the state components (say, $y_1$) is monotone increasing
along any feasible trajectory $\y(t)$, then ordering all the meshpoints
in ascending order by $x_1$ would guarantee that only one Gauss-Seidel 
iteration is needed.  
Such ordering is analogous to a simple time-marching
(from the past into the future) used with discretizations
of time-dependent PDEs (e.g., in recovering the value function
for fixed-horizon optimal control problems).
If a causal ordering is explicitly known a priori (as in the above example),
we refer to such discretizations as {\em explicitly causal}.

Explicit causality is a desirable property since it 
results in computational efficiency.  
Suppose $u(\x)$ is a viscosity solution of some boundary value problem.
Historically, two approaches have sought to capitalize on
explicit causality in solving more general static PDEs.
First, it is often possible to formulate a different time-dependent PDE 
on the same domain $\domain$ so that its viscosity solution 
$\phi$ is related to $u$ as follows:  
$u(\x) = t \; \Longleftrightarrow \; \phi(\x,t) = 0$.
The PDE for $\phi$ is then solved by explicit time-marching;
moreover, since only the zero level set of $\phi$
is relevant for approximating $u$, the Level Set Methods are applicable
e.g., see \cite{Osher_Dirichlet}, \cite{SethBook2}.  Aside from increasing 
the dimensionality of the computational domain, this approach is also
subject to additional CFL-type stability conditions, 
which often impact the efficiency substantially.  
Alternatively, in some applications (e.g., seismic imaging)
a ``paraxial'' approximation results from
assuming that all {\em optimal} trajectories must be monotone
in one of the state components (say, $y_1$) even if the same is not true
for all {\em feasible} trajectories.  This leads to a time-like
marching in $x_1$ direction (essentially solving a time-dependent PDE
on an $(n-1)$-dimensional domain); see, e.g., \cite{QianSymes1}.
This method is certainly computationally efficient, but if
the assumption on monotonicity of optimal trajectories is incorrect,
it does not converge to the solution of the original PDE.

\subsection{Efficient methods for non-explicitly causal problems.}
\label{ss:marching}

Finding the right ordering 
without increasing the dimensionality
has been the subject of much work in the last fifteen years.
This task can be challenging for discretizations which are causal,
but not explicitly causal.

For the isotropic control problems and the Eikonal 
PDE~\eqref{eq:Eikonal}, 
an equivalent of discretization~\eqref{eq:geometric_aniso} 
on a Cartesian grid is {\em monotone-causal} 
in the sense that $U_i$ cannot depend on $U_j$ 
if $U_j > U_i$.  This makes it possible to find the correct ordering
of gridpoints (ascending in $U$) at run-time, effectively
de-coupling the system of discretized equations.
This idea, the basis of Dijkstra's classic method for
shortest paths on graphs \cite{Diks}, yields the computational
complexity of $O(M \log M)$.
Tsitsiklis has introduced the first Dijkstra-like algorithm
for semi-Lagrangian discretization on a Cartesian grid in \cite{Tsitsiklis_conference}.
A Dijkstra-like Fast Marching Method was introduced by Sethian
for Eulerian upwind finite-difference discretizations of
isotropic front propagation problems \cite{SethFastMarcLeveSet}.
The method was later extended by Sethian and collaborators 
to higher-order accurate discretizations on 
grids and unstructured meshes, in $\R^n$ and on manifolds,
and to related non-Eikonal PDEs; see \cite{SethBook2},
\cite{SethVlad1}, and references therein.
For anisotropic control problems, the discretization
~\eqref{eq:geometric_aniso} generally is not causal,
but the computational stencil can be expanded dynamically  
to regain the causality.  This is the basis of
Ordered Upwind Methods \cite{SethVlad2,SethVlad3,VladThesis},
whose computational complexity is $O( \Upsilon M \log M)$,
where $\Upsilon$ measures the amount 
of anisotropy present in the problem.

The idea of alternating the directions 
of Gauss-Seidel sweeps to ``guess'' the correct mesh point ordering 
was employed by Boue and Dupuis to speed up 
the convergence in \cite{BoueDupuis}.  
For Eulerian discretizations of HJB PDEs,
the same technique was later used as a basis of Fast Sweeping Methods 
by Tsai, Cheng, Osher, and Zhao \cite{TsaiChengOsherZhao, ZhaoFSM}.
Even though a finite number of Gauss-Seidel sweeps is required
in the Eikonal case, resulting in a computational cost of $O(M)$,
the actual number of sweeps is proportional 
to the maximum number of times an optimal trajectory 
may be changing direction from quadrant to quadrant. 
A computational comparison of fast marching and fast sweeping 
approaches to Eikonal PDEs can be found in \cite{HysingTurek, GremaudKuster}.

\section{Multi-criterion optimal control.}
\label{s:multi-intro}

Unlike the classical case considered in section \ref{s:traditional},
in realistic applications there is often more than one criterion
for optimality.  
For a system controlled in $\domain \subset \R^n$, 
there may be a number of different running costs
$K_0, \ldots, K_r$ and a number of different terminal costs
$q_0, \ldots, q_r$
(all of them satisfying assumptions (A0-A3))
resulting in $(r+1)$ different definitions 
$\J_0, \ldots, \J_r$ of the total cost for a particular control.
A common practical problem is to find a control $\ba^*(\cdot)$
minimizing $\J_0(\x, \ba^*(\cdot))$ but subject to constraints 
$\J_i(\x, \ba^*(\cdot)) \leq B_i$ for all $i=1, \ldots, r$.

We will refer to a control minimizing $\J_0$ without any regard to
constraints as {\em primary-optimal}.  A control minimizing 
some $\J_j$ (for $j > 0$) without any regard to other constraints
will be called {\em $j$-optimal} (or simply {\em secondary-optimal}, 
when $j$ is clear from the context).  A control minimizing $\J_0$
subject to all of the above constraints on $\J_i$'s will be called
{\em constrained-optimal}.

For a fixed $\x \in \domain$,
we will say that a control $\ba_1(\cdot)$
is {\em j-dominated} by a control $\ba_2(\cdot)$ if
$\J_i(\x, \ba_2(\cdot)) \leq \J_i(\x, \ba_1(\cdot))$ for all $i=0, \ldots,r$
and
$\J_j(\x, \ba_2(\cdot)) < \J_j(\x, \ba_1(\cdot))$. 
We will also say that $\ba_1(\cdot)$
is {\em dominated} by $\ba_2(\cdot)$
if it is j-dominated for some $j \in \{0, \ldots, r\}$.
E.g., the constrained-optimal control $\ba^*(\cdot)$ described above
might be dominated, but will not be 0-dominated by any control.
Any control $\ba^{**}(\cdot)$ dominating $\ba^*(\cdot)$ will have
the same {\em primary cost} $\J_0$ and will also satisfy 
the same constraints; moreover, it will even spend less in at least one
of the {\em secondary costs} $\J_1, \ldots, \J_r$.

We will say that $\ba(\cdot)$ is {\em Pareto-optimal} for $\x$,
if it is not dominated by any other control.  In that case,
its vector of costs corresponds to a point on a Pareto-front 
$\mathcal{P}(\x)$ in a $\J_0, \ldots, \J_r$ space; see Figure \ref{fig:Pareto}A.

Our goal is to find an efficient numerical method for approximating 
the costs associated with Pareto-optimal controls.  
The efficiency requires solving this problem for all $\x \in \domain$
simultaneously.

We begin by showing how to compute the total cost $\J_i$ incurred 
by using a control optimal with respect to a different running cost
(section \ref{ss:secondary-along-primary}).
We then describe a recent method due to Mitchell and Sastry for
recovering a convex portion of the Pareto front by {\em scalarization}
(section \ref{ss:Mitchell-Sastry}).
Finally, in section \ref{ss:new-approach} we describe a new method
for solving fully this problem by augmenting the system state 
to include the ``budget remaining'' in each secondary cost.

\subsection{Total cost along ``otherwise optimal'' trajectories.}
\label{ss:secondary-along-primary}

We will use $u_i(\x)$ to denote the value function with respect to
$\J_i$ if all other costs are ignored.  As explained in section \ref{ss:static},
$u_i(\x)$ can be recovered as a viscosity solution of 
the PDE~\eqref{eq:HJB_general}, if we set $K=K_i$ and $q=q_i$.

Given a different value function $u$ derived for some other cost $\J$,  
we will define a restricted set of $\J$-optimal controls
$$
\mathcal{A}_{u, \blds{x}} = \left\{ \ba(\cdot) \in \mathcal{A} \, \mid \, 
\J(\x, \ba(\cdot)) = u(\x) \right\}.
$$
As explained in section \ref{ss:static}, if $\J$ satisfies the assumptions
(A0-A3), then $\mathcal{A}_{u, \blds{x}}$ will be non-empty for every $\x \in \domain$
and will contain a single element for every $\x \in \domain \backslash \Gamma$.

We will use $v_i(\x)$ to denote a $\J$-optimality-restricted
value function with respect to $\J_i$:
\begin{equation}
v_i(\x) \; = \; \inf_{\ba(\cdot) \in \mathcal{A}_{u, \blds{x}}} \J_i(\x, \ba(\cdot)).
\label{eq:value_restricted}
\end{equation}
This notation relies on a fixed choice of $\J$ and $u$, 
and we will specify $\J$ explicitly in each case to avoid ambiguity.
For the cases when $\J = \J_j$ and $u = u_j$, we will
use $v_{ij}$ instead of $v_i$.  According to this definition, we also have $v_{ii} = u_i$.

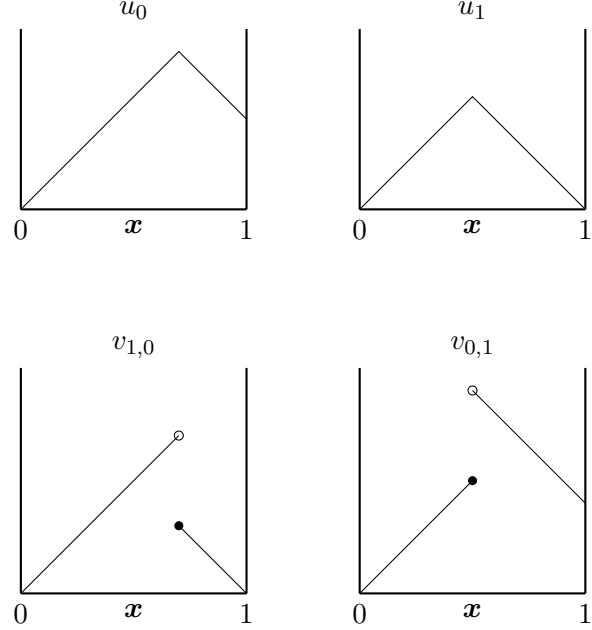
\begin{figure}[hhhh]
\center{
$$
\begin{array}{cccc}
\begin{tikzpicture}[scale=3]

	\draw[thick]			(0,0) node[below]{$0$} -- (0.5,0) node[below]{$\x$} -- (1,0) node[below]{$1$};
	\draw[thick]			(0,0) -- (0,0.8);
	\draw[thick]			(1,0) -- (1,0.8);
	\draw			(0,0) -- (0.7,0.7) -- (1,0.4);
	\draw			(0.5,0.8) node[above]{$u_0$};

\end{tikzpicture}
&&&
\begin{tikzpicture}[scale=3]

	\draw[thick]			(0,0) node[below]{$0$} -- (0.5,0) node[below]{$\x$} -- (1,0) node[below]{$1$};
	\draw[thick]			(0,0) -- (0,0.8);
	\draw[thick]			(1,0) -- (1,0.8);
	\draw			(0,0) -- (0.5,0.5) -- (1,0);
	\draw			(0.5,0.8) node[above]{$u_1$};

\end{tikzpicture}\\
\\
\\
\begin{tikzpicture}[scale=3]

	\draw[thick]			(0,0) node[below]{$0$} -- (0.5,0) node[below]{$\x$} -- (1,0) node[below]{$1$};
	\draw[thick]			(0,0) -- (0,1);
	\draw[thick]			(1,0) -- (1,1);
	\draw			(0,0) -- (0.7,0.7);
	\draw			(0.7,0.3) -- (1,0);
	\draw			(0.7,0.7) circle (.2mm);
	\fill			(0.7,0.3) circle (.2mm);
	\draw			(0.5,1) node[above]{$v_{1,0}$};

\end{tikzpicture}
&&&
\begin{tikzpicture}[scale=3]

	\draw[thick]			(0,0) node[below]{$0$} -- (0.5,0) node[below]{$\x$} -- (1,0) node[below]{$1$};
	\draw[thick]			(0,0) -- (0,1);
	\draw[thick]			(1,0) -- (1,1);
	\draw			(0,0) -- (0.5,0.5);
	\draw			(0.5,0.9) -- (1,0.4);
	\draw			(0.5,0.9) circle (.2mm);
	\fill			(0.5,0.5) circle (.2mm);
	\draw			(0.5,1) node[above]{$v_{0,1}$};

\end{tikzpicture}
\end{array}
$$
}
\caption{
{\footnotesize
Cost along ``otherwise optimal'' trajectories. A simple one-dimensional example: $\cdomain = [0,1]$,
$f = K_0 = K_1 = 1$, and $q_0(0)=q_1(0)=q_1(1)=0$, but $q_0(1)=0.4$.  
Similar discontinuities may occur even if $q_0=q_1$ (when $K_0 \neq K_1$). 
}
}
\label{fig:secondary_ex1}
\end{figure}

The following properties of Pareto fronts follow from the above definitions.
The proofs are simple and we omit them for the sake of brevity.

\begin{property}
$v_i(\x) \geq u_i(\x)$ for $\forall x \in \cdomain$ and for
any choice of $\J$ satisfying assumptions (A0-A3).
\end{property}

\begin{property}
Suppose $u(\x)$ is the value function corresponding
to some $\J$ that satisfies assumptions (A0-A3).
If $K_i$ and $q_i$ also satisfy (A0-A2),
then $v_i(\x)$ is lower semi-continuous on $\domain$
and continuous wherever $u(\x)$ is differentiable.
\end{property}

\begin{property}
Let 
$\mathcal{U}(\x) = \left\{ (J_0, \ldots, J_r) \in \R^{r+1} \, \mid \, 
J_i \geq u_j(\x), \, j=0,\ldots,r \right\}.$
Then for $\forall \x \in \cdomain$,
\begin{enumerate}
\item
$\mathcal{P}(\x) \subset \mathcal{U}(\x)$ 
and
\item
$\left( v_{i0}(\x), \ldots, v_{ir}(\x) \right) \in \mathcal{P}(\x)$
for all $i=0,\ldots,r$.
\end{enumerate}
\label{prop_Pfront_contained}
\end{property}

\begin{property}
If $r=1$, then 
$\mathcal{P}(\x) \subset [v_{00}(\x), v_{01}(\x)] \times [v_{11}(\x), v_{10}(\x)].$\\
\end{property}

Suppose $u$ is a smooth solution to the PDE~\eqref{eq:HJB_general}.
Then $\Gamma = \emptyset$ and for every $\x \in \domain$ there exists
a unique optimal/minimizing control value $\ba^* = \ba^*(\x, \grad u(x)) \in A$.
The local rate of increase of $\J_i$ along the corresponding trajectory
is $K_i(\x, \ba^*)$.  This yields a system of $(r+1)$ linear PDEs
\begin{eqnarray}
\nonumber
\nabla v_i(\x) \cdot \fB(\x, \ba^*) &=& -K_i(\x, \ba^*),
\qquad
\forall \x \in \domain;\\
v_i(\x) &=& q_i(\x),
\qquad \qquad
\forall \x \in \boundary;
\qquad i=0, \ldots, r.
\label{eq:other_optimal}
\end{eqnarray}

This system is coupled to a nonlinear equation~\eqref{eq:HJB_general},
since $\ba^*$ is generally not available a priori unless $\nabla u$ is already known.

In the Eikonal case 
(when $A= S_{n-1}$, $\fB(\x, \ba) = f(\x) \ba$ and $K(\x, \ba) = K(\x)$),
the optimal direction of motion is
$$
\ba^* = - \frac{\nabla u (\x)}{| \nabla u(\x)|} = -\nabla u(\x) \frac{f(\x)}{K(\x)},
$$
and the system~\eqref{eq:other_optimal} can be rewritten as
\begin{eqnarray}
\nonumber
\nabla v_i(\x) \cdot \nabla u(\x) &=& K_i(\x, \ba^*) K(\x) / f^2(\x),
\qquad
\forall \x \in \domain;\\
v_i(\x) &=& q_i(\x)
\qquad
\forall \x \in \boundary;
\qquad i=0, \ldots, r.
\label{eq:other_optimal_Eikonal}
\end{eqnarray}

If $u$ is not smooth, the functions $v_i$ may become discontinuous (see Figure \ref{fig:secondary_ex1})
and a generalized solution is needed to define $v_i(\x)$ at any points $\x \in \Gamma$.
Luckily, Bellman's optimality principle provides an alternative definition
to resolve this ambiguity:
$$
v_i(\x) = \lim\limits_{\tau \to 0^+} \,
\min\limits_{\ba^* \in A_{u, \blds{x}}}
\, \left\{
\tau K(\x, \ba^*) + 
v_i \left(\x + \tau \fB(\x, \ba^*) \right)
\, \right\},
$$
where $A_{u, \blds{x}} \subset A$ is the set of minimizing control values 
in~\eqref{eq:HJB_general}.  If $\x \not \in \Gamma$, then 
the $A_{u, \blds{x}}$ consists of a single element and this formulation 
is equivalent to~\eqref{eq:other_optimal}.  
Whether or not $\ba^*$ is unique, the above formula yields 
the lower semi-continuity of $v_i$. It can also be used (with a 
fixed small $\tau >0$) to derive a first-order semi-Lagrangian discretization 
of~\eqref{eq:other_optimal}.

We note that the key technical idea employed above (solving a linear equation
along the characteristics of another PDE) is well-known and useful in many 
applications.  A common version of it arises whenever there is a need
to ``propagate a constant'' from the boundary along 
the characteristics of some PDE.  
For example, this is the essential idea behind 
the ``velocity extension equation'' in \cite{SethAdalst}
and the ``escape equations'' in \cite{SethFomel}.
A slightly less general version of our equation~\eqref{eq:other_optimal_Eikonal}
(for isotropic costs/dynamics with $f=1$ and $q_i=0$) 
has also been previously used in \cite{MitchellSastry}.

\subsection{Weighted sums method and Pareto fronts.}
\label{ss:Mitchell-Sastry}

Scalarization is a popular approach where a point on Pareto-front 
is recovered by minimizing an ``aggregate objective function''.
The {\em method of weighted sums} defines that aggregate function
as a convex linear combination of the original objectives \cite{Marler}.

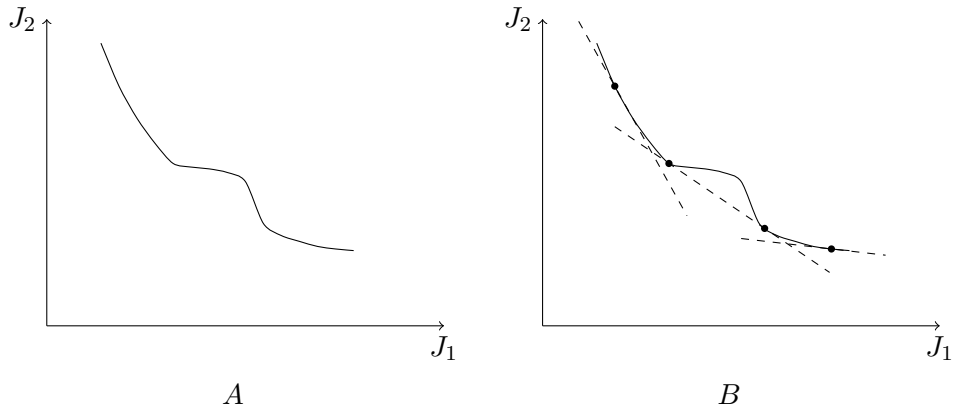
\begin{figure}[hhhh]
\center{
$$
\hspace*{-2mm}
\begin{array}{cc}
\begin{tikzpicture}[scale=2.4]

	\draw[->]		(0.3,0.1) -- (0.3,1.8)	node[left]{$J_2$};
	\draw[->]		(0.3,0.1) -- (2.5,0.1)	node[below]{$J_1$};

	\draw plot[smooth, ultra thick]	coordinates 
	{
	(.6,1.667) 	(.7, 1.4286) (.8, 1.25) (.9, 1.1111) (1,1)
	 (1.1, 0.98) (1.2, 0.97) (1.3, 0.95) (1.4, 0.9)
	 (1.5, 0.667) 
	 (1.6, 0.600) 
	 (1.7, 0.568) (1.8, 0.540) 
	 (1.9, 0.526) (2, 0.517)
};
\end{tikzpicture}
&
\begin{tikzpicture}[scale=2.4]

	\draw[->]		(0.3,0.1) -- (0.3,1.8)	node[left]{$J_2$};
	\draw[->]		(0.3,0.1) -- (2.5,0.1)	node[below]{$J_1$};

	\draw plot[smooth, ultra thick]	coordinates 
	{
	(.6,1.667) 	(.7, 1.4286) (.8, 1.25) (.9, 1.1111) (1,1)
	 (1.1, 0.98) (1.2, 0.97) (1.3, 0.95) (1.4, 0.9)
	 (1.5, 0.667) 
	 (1.6, 0.600) 
	 (1.7, 0.568) (1.8, 0.540) 
	 (1.9, 0.526) (2, 0.517)
};
	\draw[dashed] (0.7,1.204) -- (1.89, 0.396);
	\fill			(1,1) circle (.2mm);
	\fill			(1.53,0.64) circle (.2mm);
	\draw[dashed]			(1.4, 0.584) -- (2.2, 0.4915);
	\fill			(1.9, 0.526) circle (.2mm);
	\draw[dashed]	(.5, 1.788) -- (1.1, 0.71);
	\fill			(.7, 1.4286) circle (.2mm);
\end{tikzpicture}\\
A & B
\end{array}
$$
}
\caption{
{\footnotesize
Pareto Front (left) and 
its reconstruction using the ``weighted sums'' method (right).
An optimum found for each $\lambda \in \Lambda$ yields a point 
on a convex part of Pareto front and the vector $\lambda$
will be orthogonal to a support hyperplane to the front at that point
(shown by a dashed line).  In weighted sums method,
an envelope of all support hyperlplanes is used to approximate 
the Pareto front, but misses all non-convex parts of the front
\cite{Das}. 
}
}
\label{fig:Pareto}
\end{figure}

A recent method based on this approach was introduced 
by Mitchell \& Sastry for multiobjective optimal control
in the case when $f=1$ and all costs are isotropic \cite{MitchellSastry}.
Here we describe a slightly generalized version of their method.

Let
$\Lambda = \{ \lambda = (\lambda_0, \ldots, \lambda_r) \, \mid \,  
\sum\limits_{j=0}^r \lambda_j = 1  \text{ and } \lambda_i \geq 0 \text{ for all } i \}$
and suppose $\tilde{\Lambda}$ is some mesh discretizing $\Lambda$.

Given $\lambda \in \Lambda$, define 
$$
K_{\lambda}(\x,\ba) = \sum\limits_{i=0}^r \lambda_i K_i(\x, \ba);
\qquad
q_{\lambda}(\x,\ba) = \sum\limits_{i=0}^r \lambda_i q_i(\x, \ba).
$$
Solve equation~\eqref{eq:HJB_general} for $K=K_{\lambda}$ and 
$q=q_{\lambda}$.\\
%      COMMENT REMOVED
Having found $u$, solve the system~\eqref{eq:other_optimal} to obtain $v_i$'s.\\
The resulting point $\left(v_0(\x), \ldots, v_r(\x)\right)$
will belong to the Pareto front $\mathcal{P}(\x)$.

The above procedure is applied repeatedly for all meshpoints
$\lambda \in \tilde{\Lambda}$ to obtain a mesh approximating $\mathcal{P}(\x)$.
Since finding each point on the Pareto front involves solving one non-linear
and $(r+1)$ linear PDEs, any restriction of the computational domain 
$\domain \subset \R^n$ is worthwhile.  This can be often accomplished by finding
$u_1, \ldots, u_r$ first and excluding all $\x$ for which 
$u_i(\x) > B_i$ for some $i \in \{1, \ldots, r\}.$

Aside from the computational cost of the above procedure, this approach suffers 
from two usual problems associated with the method of weighted sums.
First, a uniform mesh on $\Lambda$ often results in a highly non-uniform
mesh on $\mathcal{P}(\x)$.  Secondly, the weighted sum method can approximate
only a convex part of the Pareto front \cite{Das}; see Figure \ref{fig:Pareto}B.  
This may result in selecting suboptimal trajectories.  
Mitchell \& Sastry acknowledge that, for non-convex fronts,
``this method may fail to detect a feasible path even if one exists'', but they
report that ``non-convexity has not been a problem''
in their numerical experiments. They further suggest that 
``neighboring values of $\lambda$ can be used to bound the error in
the convex approximation'' of non-convex portions of the front.
We believe that the latter procedure can be very inaccurate, especially
since the Pareto front is frequently discontinuous for multiobjective 
optimal control problems.  See Figure \ref{fig:pareto} for an example of 
non-convexities actually present in the test-problem introduced 
in \cite{MitchellSastry}.
  
In addition, we note that recovering the entire Pareto front for each 
$\x \in \domain$ is excessive and unnecessary when the real goal is to solve
the problem for a fixed list of constraints $B_1, \ldots, B_r$.
If instead of using a mesh $\tilde{\Lambda}$, we attempt changing 
$\lambda$ adaptively, it is not generally clear in what direction should
$\lambda$ be perturbed in order to satisfy the constraints.

Despite the above limitations, this technique can be useful in many applications,
whenever there is a need to produce at least some of the Pareto-optimal 
trajectories efficiently.  E.g., a similar scalarization approach 
has been independently used by Kim and Hespanha in path planning 
(to minimize the risk of detection/interception) for groups of UAVs \cite{KimHespanha}.

\subsection{An augmented PDE on an expanded state space.}
\label{ss:new-approach}

We propose an alternative approach based on augmenting 
the system-state space to reflect the budget remaining in each 
of the secondary costs.  This is a generalization of the idea
classically used to recast a Bolza problem as a problem with 
zero running cost by adding an extra dimension to the state space \cite{BressanPiccoli}.

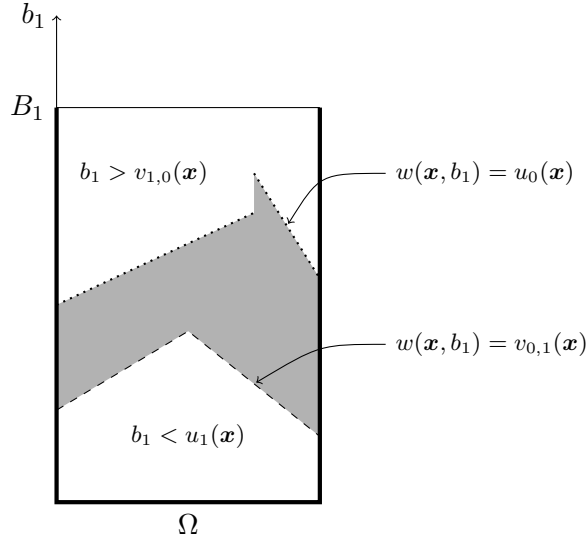
\begin{figure}[hhhh]
\center{
\begin{tikzpicture}[scale=1.75]

	\fill[gray!60!white]	(0,1.5) -- (1.5,2.2) -- (1.5,2.5) -- (2,1.7) --
							(2,0.5) -- (1,1.3) -- (0,0.7) -- cycle;				

	\draw[->]		(0,0) -- (0,3.7)	node[left]{$b_1$};
%      COMMENT REMOVED
	\draw[ultra thick]		(0,3) node[left]{$B_1$} -- (0,0) -- 
							(1,0) node[below]{$\domain$} --
							(2,0) -- (2,3);
	\draw			(0,3) -- (2,3);
	\draw[dashed]			(0,0.7) -- (1,1.3) -- (2,0.5);
	\draw[dotted,thick]		(0,1.5) -- (1.5,2.2);
	\draw[dotted,thick]		(1.5,2.5) -- (2,1.7);

	\draw[font=\footnotesize]	(1,0.5) node {$b_1 < u_1(\x)$};
	\draw[font=\footnotesize]	(0.1,2.5) node[right] {$b_1 > v_{1,0}(\x)$};

	\draw[font=\footnotesize, -to, shorten >= 1pt]	(2.5,2.5) node[right] {$w(\x,b_1) = u_0(\x)$} .. controls(2,2.5) .. (1.75,2.1);
	\draw[font=\footnotesize, -to, shorten >= 1pt]	(2.5,1.2) node[right] {$w(\x,b_1) = v_{0,1}(\x)$} .. controls(2,1.2) .. (1.5,0.9);

\end{tikzpicture}
}
\caption{
{\footnotesize
A sketch of the domain $\cdomain_e$ for the case $r=1$.  
The thick line shows the inflow boundary $\inBoundary$.
The dashed line is the graph of $u_1(\x)$ 
(the minimum $b_1$ needed to leave $\domain$ starting from $\x$).
$\fBoundary$ is the portion of $\inBoundary$ above the dashed line.
The dotted line shows the points where $b_1 = v_{1,0}(\x)$.
At that level and higher, the primary-optimal controls are feasible
and $w(\x,b_1) = u_0(\x)$.  Above the dotted line, the constraint is slack.
The PDE \eqref{eq:multiobject_PDE}
has to be solved numerically in the shaded region.
}
}
\label{fig:_domain}
\end{figure}

Suppose $b_i \in [0, B_i]$ is the ``budget'' remaining in the secondary cost $\J_i$.
We define an extended state variable $\xhat = (\x, b_1, \ldots, b_r)$
and the extended state space 
$\domain_e = \domain \times (0, B_1) \times \ldots \times (0, B_r)$.
The {\em outflow boundary} and the {\em inflow boundary} of this domain are
$$
\outBoundary \; = \; \left\{
\xhat \in \partial \cdomain_e \, \mid \,
\x \in \domain, \, \text{ and } b_i \in (0, B_i]
\text{ for } i=1, \ldots, r\right\}
\quad
\text{and}
\quad
\inBoundary \; = \; \boundary_e \backslash \outBoundary.
$$
We note that $\outBoundary$ is the part of the boundary
where at least one of the budgets is at the maximal level 
($\exists j$ such that $b_j = B_j$).
For the case $r=1$, $\inBoundary$ coincides with 
the so called {\em parabolic boundary} \cite{Evans}; see Figure \ref{fig:_domain}.
We will also refer to a {\em feasible subset} of $\inBoundary$:
$$
\fBoundary = \left\{
\xhat \in \inBoundary \, \mid \,
\x \in \boundary, \text{ and } b_i \geq q_i(\x) 
\text{ for } i = 1, \ldots, r
\right\}.
$$
The extended state at the time $t$ will be denoted as 
$\yhat(t) = (\y(t), \beta_1(t), \ldots, \beta_r(t)) \in \bar{\domain}_e,$
and the extended dynamics is now defined by
\begin{eqnarray}
\nonumber
\y^{\prime}(t) &=& \fB(\y(t), \ba(t));\\
\beta_i^{\prime}(t) &=& -K_i(\y(t), \ba(t)), \qquad i=1,\ldots,r;\\
\yhat (0) &=& \xhat = (\x, b_1, \ldots, b_r) \in \domain_e.
\label{eq:extended_dynamics}
\end{eqnarray}
As before, the exit time corresponding to a particular control is defined by~\eqref{eq:exit-time},
but we will use $T=T_{\xhat, \ba}$ to simplify the notation.

The total cost of this control is defined as
\begin{equation}
\Jhat(\xhat, \ba(\cdot)) = 
\begin{cases}
\int_0^{T} K_0 (\y(t), \ba(t)) \, dt 
\, + \, q_0(\y(T)) &
\text{ if } \yhat(T) \in \fBoundary,\\
+\infty & \text{ otherwise.}
\end{cases}
\label{eq:multiobject_control_cost}
\end{equation}
The (possibly infinite) value function of the new control problem is defined, as usual:
$w(\xhat) = w(\x, b_1, \ldots, b_r) = \inf \Jhat(\xhat, \ba(\cdot)).$ 

If the value function is smooth, the standard argument (based on Taylor-expanding
$w\left(\yhat(\tau)\right)$ in the optimality principle), shows that $w$ satisfies the PDE
\begin{equation}
\min\limits_{\ba \in A}
\left\{ 
K_0 (\x, \ba) + \nabla_{\blds{x}} w \cdot \fB ( \x, \ba )
\, - \,
\sum\limits_{i=1}^r 
K_i(\x, \ba) \frac{\partial w}{\partial b_i} 
\right\}
\, = \, 0.
\label{eq:multiobject_PDE}
\end{equation}
with the boundary conditions 
\begin{equation}
w(\xhat) = 
\begin{cases}
q_0(\x) & \text{ on $\fBoundary$;}\\
+ \infty & \text{ on $\inBoundary \backslash \fBoundary.$}
\end{cases}
\label{eq:multiobject_boundary_cond}
\end{equation}
However, $w$ can be not only non-smooth, but also discontinuous
inside $\domain_e$
since the boundary data is discontinuous and no local controllability
(in the directions $b_1, \ldots, b_r$)
is present; see Remark \ref{rem:boundary_continuity}.
%      COMMENT REMOVED
Nevertheless, $w$ can still be interpreted as a unique
{\em discontinuous viscosity solution} \cite{Soravia_IntegralConstraint}
of equation~\eqref{eq:multiobject_PDE},
despite the fact that Soner's ``inward pointing'' condition
is clearly violated on $\outBoundary$ \cite{MottaRampazzo}.

As in the single-objective case, 
if the minimization in $\ba$ can be performed analytically, 
the augmented PDE can be rewritten in a simpler form.  
For example, in the fully isotropic case,
\begin{equation}
K_0 (\x) \, - \, |\nabla_{\blds{x}} w| f(\x)
\, - \,
\sum\limits_{i=1}^r 
K_i(\x) \frac{\partial w}{\partial b_i} 
\; = \; 0,
\label{eq:multiobject_Eikonal}
\end{equation}
which in the case $r=0$ reduces to the usual Eikonal equation~\eqref{eq:Eikonal}.

The following properties of $w$ follow from the above definitions.
The proofs are simple or standard and we omit most of them
(except for the last three) for the sake of brevity.

\begin{property}
The value function 
$w: \cdomain_e \mapsto \left(\R \cup \{+\infty\} \right)$ 
is lower-semicontinuous; 
see \cite{Soravia_IntegralConstraint} and references therein.
\end{property}

\begin{property}
The value function $w$ is monotone non-increasing  in each $b_i$.
As a result, if $\bb \geq \bc$ componentwise, then $w(\x,\bb) \leq w(\x, \bc).$
\end{property}

\begin{property}
If for some $i$, $b_i < u_i(\x)$, then $w(\x,\bb) = +\infty.$
\label{prop:min_bi_needed}
\end{property}

\begin{property}
$w(\x, \bb) \geq u_0(\x) = 
w \left(
\x, v_{10}(\x), \ldots, v_{r0}(\x) 
\right).$
\label{prop:w_lower_bound}
\end{property}

\begin{property}
$w$ is Lipschitz continuous along every characteristic.
\label{prop:w_Lip_on_chars}
\end{property}

\begin{property}
The characteristics of PDE ~\eqref{eq:multiobject_PDE} have the following properties:
\begin{enumerate}
\item
All characteristics are monotone-increasing in all $b_i$'s.
\item
Projections of characteristics on $\domain$ yield
constrained-optimal trajectories.
\item
Any characteristic starting on $\fBoundary$ leads into $\cdomain_e$.
\item
Any characteristic starting on $\outBoundary$ leads out of $\cdomain_e$.
\end{enumerate}
\label{prop:causal_chars}
\end{property}

We will say that a control $\ba(\cdot)$ is {\em feasible} 
for $\xhat = (\x, \bb)$ if $\J_i(\x, \ba(\cdot)) \leq b_i$ for all $i=1, \ldots,r$.
We will say that $\xhat \in \cdomain_e$ is a {\em feasible point}
if it has at least one feasible control (i.e., if $w(\xhat) < \infty$).
We will say that a point $\xhat$ is {\em $i$-tightly-constrained} if
for every constrained-optimal control $\ba^*(\cdot)$ we 
have $\J_i(\x, \ba^*(\cdot)) = b_i$.  Otherwise, we will
call $\xhat$ {\em $i$-slack}.  We will also say that $\xhat$
is {\em totally slack} if there exists a constrained-optimal control 
$\ba^*(\cdot)$ such that $\J_i(\x, \ba^*(\cdot)) < b_i$
for all $i=1, \ldots, r.$

\begin{property}
The point $\left(w(\x, \bb), \bb \right)$ lies
on the Pareto front $\mathcal{P}(\x)$ in $\R^{r+1}$
if and only if $(\x, \bb)$ is $i$-tightly-constrained for all $i=1, \ldots, r.$\\
Moreover, if $(b_0, \ldots, b_r) \in \mathcal{P}(\x)$ and $b_i \leq B_i$
for all $i=1, \ldots, r$, then $w(\x, b_1, \ldots, b_r) = b_0$.
\label{prop:w_and_Pareto}
\end{property}
\begin{proof}
Suppose $\ba^*(\cdot)$ is the constrained-optimal control for $\xhat=(\x, \bb).$
If $\xhat$ is $i$-slack, then $\J_i(\x, \ba^*(\cdot)) < b_i$ and this does
not contribute a point on Pareto front.

If $\xhat$ is $i$-tightly-constrained for all $i=1, \ldots, r,$
then $\ba^*(\cdot)$ cannot be dominated by any control.
(If it were 0-dominated, that would contradict its constrained-optimality.
If it were $i$-dominated for $i>0$, that would contradict the fact that
$\xhat$ is $i$-tightly-constrained.)

Finally, suppose $\ba(\cdot)$ is a control that realizes the
cost vector $(b_0, \ldots, b_r)  \in \mathcal{P}(\x)$ with $b_i \leq B_i$
for all $i=1, \ldots, r$.  By the definition of $w$, we have
$w(\x,b_1,\ldots,b_r) \leq \J_0(\x, \ba(\cdot)) = b_0$.
If we had $w(\x,b_1,\ldots,b_r) <  b_0$ this would imply the existence
of a control $0$-dominating $\ba(\cdot)$, which would contradict its Pareto-optimality.
So,  $w(\x,b_1,\ldots,b_r) =  b_0$.
\end{proof}

\begin{property}
If the point $(\x_1,\bb)$ is totally slack, 
then $w$ is locally Lipschitz-continuous in its first argument.
\label{prop:w_slack_Lip}
\end{property}
\begin{proof}
First, we note that there exists some open neighborhood of $\x_1$
such that $(\x_2, \bb)$ is totally slack for every $\x_2$ from 
that neighborhood.  Otherwise there would exist a sequence of 
$i$-tightly-constrained $\x_j$'s converging to $\x_1$, 
which can be used to show that $(\x_1,\bb)$ is $i$-tightly-constrained as well.

If $(\x_1,\bb)$ is totally slack, then, starting from any point $\x_2$
close to $\x_1$, the $i$-th cost of reaching $\x_1$ is bounded above by 
$ |\x_1 - \x_2| \max (K_i / |\fB|)$,
where the maximum is taken over all $\x$ and all $\ba$ such that
$\fB(\x,\ba) / K_i(\x,\ba) \in \partial V_i(\x)$.
By the assumption (A3), $V_i$ contains the origin in its interior;
so, there exists a constant $L$ such that the $i-$th cost of reaching
$\x_1$ from $\x_2$ is bounded by $L |\x_1 - \x_2|$ for all $i$.
Suppose we travel from $(\x_2, \bb)$ to $(\x_1, \btilde)$, where
$\btilde \leq \bb$.
If $\x_2$ is sufficiently close to $\x_1$, then $(\x_1, \btilde)$
is still totally slack, and any constrained-optimal control 
for $(\x_1, \bb)$ will be still feasible for $(\x_1, \btilde)$.
Thus, $w(\x_2, \bb) \leq L |\x_1 - \x_2| + w(\x_1, \bb)$.
To complete the proof, we repeat the argument
switching the roles of $\x_1$ and $\x_2$.
\end{proof}

\begin{property}
If $\xhat = (\x,\bb)$ is totally slack, $\ba^*(\cdot)$ is a constrained-optimal control for 
$\xhat$ and $\y^*(t)$ is the corresponding trajectory in $\cdomain$,
then $\y^*(t)$ is also a characteristic of problem~\eqref{eq:HJB_general} 
with $K=K_0$ and $q=q_0$.
\end{property}
\begin{proof}
Briefly, if $(\x,\bb)$ is totally slack, then any sufficiently small perturbation 
of $\ba^*(\cdot)$ will also be feasible.  Since $\ba^*(\cdot)$ is constrained-optimal,
the function $\y^*(t)$ is a local minimizer of $\J_0$ and will be a solution of 
the Euler-Lagrange equation in $\R^n$ (see, e.g., \cite{Evans}).  
By Pontryagin's maximum principle, it will also be a characteristic of 
the corresponding HJ PDE on $\cdomain$.  This is an interesting fact,
since the characteristics of ~\eqref{eq:HJB_general} yield locally 
primary-optimal (unconstrained) trajectories.
\end{proof}

Property \ref{prop:causal_chars}.1 is the basis for
explicitly causal discretizations of the augmented PDE,
which enables an efficient numerical treatment (by ``marching'' in any direction $b_i$).
Properties 
\ref{prop:min_bi_needed} and
\ref{prop:w_lower_bound}
can be used to reduce the computational domain,
as shown in Figure \ref{fig:_domain} and further explained 
in section \ref{ss:reduced_domain}.  Property \ref{prop:w_and_Pareto}
can be used to extract the relevant part of the Pareto front 
from the values of $w$.

\subsection{Numerical method for the augmented PDE}
\label{ss:new-numerics}

We consider a Cartesian grid $\hat{X}$ on $\cdomain_e$.
For simplicity, we will assume the same grid spacing $h$ in
all spatial dimensions and grid spacing 
$\Delta b_1, \ldots, \Delta b_r$ in each of the constraint/secondary 
cost directions.  Let $\hat{h} = \max \{h, \, \max_{i} \Delta b_i \}.$
Our goal is to construct an approximate
solution $W$ converging to the lower semi-continuous value function
$w$ as $\hat{h} \to 0.$ 

If the minimization in $\ba$ can be performed analytically
(e.g., in the fully isotropic case of equation~\eqref{eq:multiobject_Eikonal}), 
a natural Eulerian framework scheme
may be obtained by using upwind finite differences to approximate the partial
derivatives of $w$.  However, this approach, even when feasible, 
would be subject to CFL-type stability conditions,
which would potentially have a significant impact on the computational cost of the scheme.
(A simple example to illustrate this:  suppose the $\fB$ and $K_0$ are isotropic,
$r=1$, $q_1 = 0$ on $\boundary$, and $K_1 = 1$,
making $\J_1$ be the total time to $\boundary$ along a given trajectory.
If upwind finite differences are used in a given $b_1$-slice to approximate 
$\nabla_{\blds{x}} w$ and the forward-difference is used to approximate $\frac{\partial w}{\partial b_1}$,
this results in a scheme suitable for explicit forward ``time''-marching in $b_1$-direction,
but not surprisingly requires a standard CFL condition $\max_{\x} f(\x) \leq h / \Delta b_1$
for stability.)

Instead, we chose to implement a semi-Lagrangian discretization of
the augmented PDE~\eqref{eq:multiobject_PDE}.  
In addition to improved stability properties,
the resulting scheme is also easy to extend to unstructured meshes 
in $\cdomain_e$ and is applicable when the minimization in $\ba$ cannot be
handled analytically (which is often the case for control-theoretic problems).
Our discretization is a variant of ~\eqref{eq:Falcone-style}.
Given a point $\xhat = (\x, \bb) = (x_1, \ldots, x_n, b_1, \ldots, b_r)$, 
we define a new system state $(\xtilde, \btilde)$ as follows:
\begin{equation}
\xtilde = \x + \tau_{\blds{a}} \fB(\x, \ba); \qquad 
\tilde{b}_i = b_i - \tau_{\blds{a}} K_i(\x, \ba), 
\; \text{ for } i=1, \ldots, r.
\label{eq:newstate_after_one_step}
\end{equation}
Here the obvious dependence of the new state on the choice of control $\ba$
is omitted for the sake of notational simplicity.

\begin{equation}
W(\x, b_1, \ldots, b_r) \; = \;  
\min\limits_{\ba \in A}
\left\{
\tau_{\blds{a}} K_0(\x, \ba) \; + \; W(\xtilde, \btilde)
\right\},
\qquad \forall \text{ gridpoints } (\x, b_1, \ldots, b_r) \not \in \inBoundary.
\label{eq:discr_augmented}
\end{equation} 

To ensure that the discretized equations allow efficient marching in the direction
$b_1$, it is sufficient to take 
\begin{equation}
\tau_{\blds{a}}  \; = \;  \Delta {b_1} / K_1(\x, \ba),
\label{eq:tau}
\end{equation}
which guarantees that $\tilde{b}_1 = b_1 - \Delta b_1$ 
for any choice of $\ba$.  This means that the new position 
$(\xtilde, \btilde)$ will be in the previous ``$b_1$-slice'', in which
the values of $W$ were already computed.

Of course, $(\xtilde, \btilde)$ will usually not be a gridpoint and 
$W(\xtilde, \btilde)$ has to be approximated using the values of $W$ at nearby 
gridpoints.  Our implementation uses a standard 
tensor product of linear interpolations in all $\x$ and $\bb$ variables; 
see, e.g., \cite{CarliniFalconeFerretti}.  

For example, if $n=2$ and $r=1$, this yields a bilinear interpolation. 
Suppose $\xtilde$ lies in a cell with vertices
$\x_1, \ldots, \x_4$ (enumerated clockwise, starting from the lower left corner).  
Let $(\gamma_1, \gamma_2) =  \left(\xtilde - \x_1 \right) / h.$  Then
\begin{equation}
\begin{array}{lll}
W(\xtilde, \tilde{b}_1) &=& 
\gamma_1 \left( \gamma_2 W(\x_3, \tilde{b}_1) + (1-\gamma_2) W(\x_2, \tilde{b}_1) \right)\\
&+&
(1-\gamma_1) \left( \gamma_2 W(\x_4, \tilde{b}_1) + (1-\gamma_2) W(\x_1, \tilde{b}_1) \right).
\end{array}
\label{eq:bilinear_interp}
\end{equation}

\begin{remark}

If the value function is smooth on the cell, 
the resulting interpolation error is $O(\htilde^2)$,
where $\htilde$ is the size of that cell  
(i.e., $\htilde =  \max \{h, \, \max\limits_{i>1} \Delta b_i \}$).  
However, in the worst case, $w$ may be discontinuous,
resulting in a $O(1)$ interpolation error. 
(We note that the property \ref{prop:w_Lip_on_chars} ensures
the Lipschitz continuity of $w$ on the characteristic itself,
but the interpolation cell containing $\xtilde$ might still include a discontinuity.)

The convergence of semi-Lagrangian
schemes to discontinuous viscosity solutions is more subtle.
In \cite{BardiFalconeSoravia1,BardiFalconeSoravia2}
Bardi, Falcone and Soravia have proven that, on any compact subset where $w$
is continuous, the semi-Lagrangian approximation converges to $w$ uniformly,
provided $\htilde / \tau \to 0$ as $\tau \rightarrow 0$.
The resulting schemes were successfully used to approximate 
discontinuous viscosity solutions both in the context of infinite-horizon 
optimal control and in differential games.
For viscosity solutions continuous on the entire domain, 
a different proof \cite{BardiFalconeSoravia2} yields 
the convergence (and the convergence rate estimates) under
a milder assumption that $\htilde/\tau$ remains bounded as $\tau \rightarrow 0$.

Since we are setting $\tau$ for each $\x$ and $\ba$ separately,
the above condition for convergence to discontinuous solutions
becomes: $\htilde / (\min \tau_{\blds{a}}) \to 0$ as $(\max \tau_{\blds{a}}) \rightarrow 0$.
If $\tau_{\blds{a}}$ is selected according to \eqref{eq:tau},
this condition requires that the grid refinement should be performed
in such a way that $\htilde = o(\Delta b_1)$.  
(An alternative approach would be to
keep $\htilde = O(\Delta b_1)$, but pick $\tau_{\blds{a}}$ so that 
$ \tau_{\blds{a}} K_1(\x, \ba) = b_1 - \tilde{b}_1 = m \Delta b_1$, 
where the integer $m \to \infty$ and $m \Delta b_1 \to 0$
as $\hat{h} \to 0$.)
However, for the examples considered in this paper, 
the numerical evidence suggests convergence
even with $m=1$ (i.e., as prescribed by formula \eqref{eq:tau})
and with $\htilde = O(\Delta b_1)$.  The discontinuities 
are smeared in a narrow band with the width of that smearing
proportional to $\Delta b_1$;
e.g., see the convergence study 
for a simple example in section \ref{ss:convergence}.
We note that any higher order accurate semi-Lagrangian scheme
would need to use a larger stencil to interpolate $W(\xtilde, \btilde)$,
and an ENO or WENO reconstruction would be needed to handle 
the discontinuities.
\label{rem:interpolation_error}
\end{remark}

Even though it is not necessary in principle,
our current implementation assumes the geometric dynamics 
(defined in section \ref{ss:static}).
The minimization in formula~\eqref{eq:discr_augmented} is performed numerically
using the standard ``golden section search'' algorithm.
Once $W$ is computed, optimal controls and trajectories
are recovered by following the characteristics of PDE~\eqref{eq:multiobject_PDE} numerically.

\subsection{Reducing the computational domain \& initialization}
\label{ss:reduced_domain}

Given the high dimensionality of $\cdomain_e$, any reduction
of the domain is likely to result in substantial
savings in the computational cost.  
Two such reductions are obviously worthwhile; see Figure \ref{fig:_domain}.
In both cases we recover a surface consisting of special characteristics
of~\eqref{eq:multiobject_PDE} by efficiently solving other PDEs 
on lower-dimensional domains.
We note that a similar approach was previously proposed in \cite{VladTimeD}
for time-optimal control in the case of time-dependent dynamics.

First, by property \ref{prop:w_lower_bound}, if $w(\x,\bb) = u_0(\x)$, 
a primary-optimal control is already feasible and further increase in secondary
budgets will not yield any reduction of $w$.  
The $n$-dimensional surface on which this happens can be found a priori
(by numerically approximating functions $v_{i0}$ for $i=1, \ldots r$).

Secondly, only a subset of $\cdomain_e$ is feasible:
if the initial budget-vector $\bb$ is insufficient to reach $\boundary$
starting from $\x$, then $w(\x, \bb) = +\infty$.
We formally define the Minimal Feasible Level (MFL) as a graph of the function 
$$ 
w_1(\x, b_2, \ldots, b_r) = \min \left\{
b_1 \, \mid \, w(\x, b_1, b_2, \ldots, b_r) < +\infty \right\}.
$$
As described below, the MFL can be efficiently recovered for any $r$ by solving 
a sequence of PDEs on lower-dimensional domains.  
However, for $r=1$, this task is particularly simple, since in that case the MFL coincides
with the graph of $u_1(\x)$, and the latter can be approximated by solving 
the PDE~\eqref{eq:HJB_general} numerically on $\domain$.
As explained in \cite{Soravia_IntegralConstraint}, 
the augmented PDE~\eqref{eq:multiobject_PDE} then has to be solved 
on the epigraph of $u_1$.
Thus, using the results in section \ref{ss:secondary-along-primary},
the value\footnote{
We emphasize that the MFL does not really provide any 
additional boundary conditions for the equation~\eqref{eq:multiobject_PDE}
(since this surface consists of characteristics of that PDE).
But for a semi-Lagrangian discretization, the numerical values on 
the MFL are needed -- for an $\xhat$ just above that surface, 
the optimal $\xtilde$ might well lie in a cell 
one of whose corners is on the MFL.
We note that any $\xtilde$ falling below the MFL is obviously 
non-optimal; when such situation arises during the minimization
in equation~\eqref{eq:discr_augmented}, we use $W(\xtilde) = +\infty$.} 
of $w$ on the MFL is provided by $v_{01}(\x)$.  

\begin{remark}
To represent MFL on the grid $\hat{X}$,
our current implementation uses the smallest
integer $i$ such that $u_1(\x) \leq i \Delta b_1$
and then initializes $W(\x, i \Delta b_1) = v_{01}(\x)$.
This procedure is conservative (in the sense that it overestimates the minimal
necessary budget), but it also introduces additional errors of order $O(\Delta b_1)$;
see the convergence study in section \ref{ss:convergence}.
A better implementation can be built by locally altering the interpolation
procedure near the MFL.
\label{rem:MFL}
\end{remark}

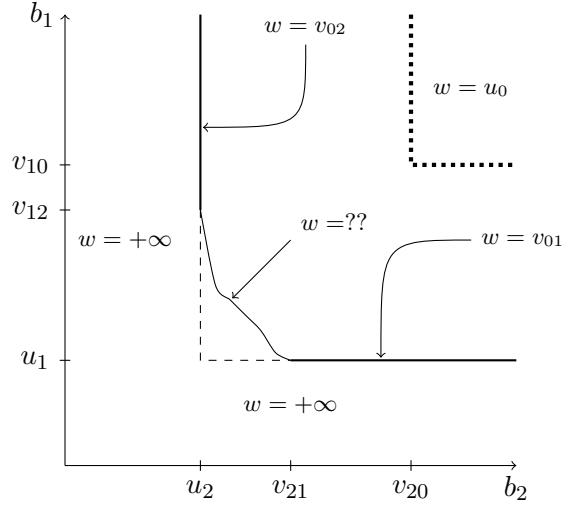
\begin{figure}[hhhh]
\center{
\begin{tikzpicture}[scale=2]

	\draw[->]		(0,0) -- (3,0)	node[below]{$b_2$};
	\draw[->]		(0,0) -- (0,3)	node[left]{$b_1$};

	\draw[thick]	(1.5,0.7) -- (3,0.7); 
	\draw[thick]	(0.9,1.7) -- (0.9,3); 

	\draw plot[smooth, thick]	coordinates 
	{ (0.9, 1.7) (1, 1.2) (1.1, 1.1) (1.2, 1) (1.3, 0.9) (1.4, 0.75) (1.5,0.7)	};

	\draw[dashed]	(0.9,1.7) -- (0.9,0.7) -- (1.5, 0.7);

	\draw[font=\footnotesize]	(1.5,0.4) node {$w = + \infty$};
	\draw[font=\footnotesize]	(0.4,1.5) node {$w = + \infty$};

	\draw		(1pt, 0.7cm) -- (-1pt, 0.7cm) node[left] {$u_1$};
	\draw		(1pt, 1.7cm) -- (-1pt, 1.7cm) node[left] {$v_{12}$};
	\draw		(1pt, 2cm) -- (-1pt, 2cm) node[left] {$v_{10}$};
	\draw		(0.9cm, 1pt) -- (0.9cm, -1pt) node[below] {$u_2$};
	\draw		(1.5cm, 1pt) -- (1.5cm, -1pt) node[below] {$v_{21}$};
	\draw		(2.3cm, 1pt) -- (2.3cm, -1pt) node[below] {$v_{20}$};

	\draw[font=\footnotesize, -to, shorten >= 1pt]	(1.6,2.8) node[above] {$w = v_{02}$} .. controls(1.6,2.25) .. (0.9,2.25);
	\draw[font=\footnotesize, -to, shorten >= 1pt]	(2.7,1.5) node[right] {$w = v_{01}$} .. controls(2.1,1.5) .. (2.1,0.7);
	\draw[font=\footnotesize, -to, shorten >= 1pt]	(1.5,1.5) node[above right] {$w = ??$} -- (1.1,1.1);

	\draw[ultra thick, dotted]	(2.3,3) -- (2.3,2) -- (3, 2);
	\draw[font=\footnotesize]	(2.7,2.5) node {$w = u_0$};

\end{tikzpicture}
}
\caption{
{\footnotesize
%      COMMENT REMOVED
Two secondary costs diagram for a single point $\x \in \domain$.
The boundary of the feasible set ($\x$'s MFL) is shown by a solid line.  
To solve the PDE using an efficient semi-Lagrangian scheme, 
it is necessary to find the entire MFL 
(including its curved and possibly non-convex part) 
and pre-initialize $w$ on that boundary.
The dotted line shows the boundary of the set 
where the unconstrained (primary-optimal) controls are feasible.
}
}
\label{fig:_2_secondary_costs}
\end{figure}

For $r>1$ the approximation of the MFL is more subtle.
In \cite{Soravia_IntegralConstraint}, Soravia suggested solving an augmented
PDE on the set $\{(\x, \bb) \in \cdomain_e \mid b_i \geq u_i(\x), \; i=1,\ldots r \}$.
However, we note that $w$ need not be finite everywhere on that set.
This is already evident for $r=2$.
Indeed, a control optimal according to $\J_1$ usually incurs a higher
$\J_2$-cost than the control optimal according to the latter (and vice versa). 
As a result, $w \left( \x, u_1(\x), u_2(\x) \right)$ is usually $+\infty.$
This situation is illustrated in Figure \ref{fig:_2_secondary_costs}
(i.e., the region between the dashed and solid lines is not feasible).
If $b_2 \geq v_{21}(\x)$, then $w(\x, u_1(\x), b_2) = v_{01}(\x)$; 
this allows to recover a part of the MFL quite easily 
(the thick solid lines in Figure \ref{fig:_2_secondary_costs}).
However, recovering the rest of the MFL
requires more computational effort.

Suppose $r=2$ and $w_1(\x, b_2)$ is the minimum budget $b_1$ 
sufficient to reach $\boundary$ from $\x$ with $\J_2 \leq b_2$. 
It is clear that the PDE for $w$ should be solved on the epigraph of $w_1$.
Since $w_1$ optimizes $\J_1$ subject to a constraint on $\J_2$,
we note that $w_1$ plays the same role for the pair of costs $(\J_1, \J_2)$
as $w$ played for the triple $(\J_0, \J_1, \J_2)$. 
Thus, the same argument used in section \ref{ss:new-approach}
shows that $w_1$ can be recovered as
a viscosity solution of a PDE similar to~\eqref{eq:multiobject_PDE},
but posed on $\cdomain \times [0, B_2]$;  i.e.,
\begin{equation}
\min\limits_{\ba \in A}
\left\{ 
K_1 (\x, \ba) + \nabla_{\blds{x}} w_1 \cdot \fB ( \x, \ba )
\, - \,
K_2(\x, \ba) \frac{\partial w_1}{\partial b_2} 
\right\}
\, = \, 0.
\label{eq:2_cost_init}
\end{equation}
To solve the latter PDE efficiently on an $(n+1)$-dimensional domain,
we first need to approximate its own $n$-dimensional ``$MFL_1$''.  Luckily, this 
is equivalent to the one-secondary-cost problem already
considered above.  Thus, this new $MFL_1$ is the graph of $u_2(\x)$
and the value of $w_1$ on the $MFL_1$ is provided by $v_{12}(\x)$.
While computing $w_1$, we can also integrate $K_0$ along $w_1$'s 
optimal trajectories (as described in section \ref{ss:secondary-along-primary})
to find the values of $w$ on its MFL
(i.e., on the graph of $w_1$ in  $\cdomain \times [0, B_1] \times [0, B_2]$).
This is the approach used to compute the example 
in section \ref{ss:two_secondary}.

For the general case ($r>2$), the same procedure can be applied recursively.
Let $w_i(\x, b_{i+1}, \ldots, b_r)$ be the minimum budget $b_i$
sufficient to reach $\boundary$ from $\x$ with $\J_j \leq b_j$
for all $j >i$.  
According to this definition, $w_r(\x) = u_r(\x)$ and,
for $ \forall i < r$, $w_i$ is finite only on the epigraph of $w_{i+1}$.
Focusing only on the costs $\J_i, \ldots \J_r$,
the argument of section \ref{ss:new-approach} shows that
$w_i$ is the (unique lower semi-continuous)
viscosity solution of
\begin{equation}
\min\limits_{\ba \in A}
\left\{ 
K_i (\x, \ba) + \nabla_{\blds{x}} w_i \cdot \fB ( \x, \ba )
\, - \,
\sum\limits_{j=i+1}^r 
K_j(\x, \ba) \frac{\partial w_i}{\partial b_j} 
\right\}
\, = \, 0
\label{eq:multiobject_recursive}
\end{equation}
on $\cdomain \times [0, B_{i+1}] \times \ldots \times [0, B_r]$.
In this $(n+r-i)$-dimensional domain,
the $MFL$ for $w_i$ (i.e., the $MFL_i$) is provided
by a graph of $w_{i+1}$, and  the values of $w_i$ on $MFL_i$
can be computed by integrating $K_i$ along the optimal
trajectories of $w_{i+1}$.
Finally, the graph of $w_1$ is the MFL for $w=w_0$,
and the latter solves the PDE~\eqref{eq:multiobject_PDE}. 

\begin{remark}
Both techniques presented in this section 
allow a significant reduction of the computational domain.
The net effect of the first technique (based on considering
only those gridpoints where $w(\x, \bb) \geq u_0(\x)$) 
depends on how large the budget bounds $B_i$s are compared
to the values of $v_{i0}(\x)$.  
If $r=1$, the percent of excluded gridpoints can be found
by averaging the ratio 
$\left[ \max \left\{ B_1 - v_{10}(\x), \, 0 \right\} \, / \, B_1 \right]$ 
over $\domain$.

The net effect of the second technique (based on finding the MFL)
depends on the degree of ``primary-non-optimality'' of
``secondary optimal'' trajectories.
If $r=1$, the percent of gridpoints additionally excluded 
(after the MFL is computed) can be found by averaging the ratio
$\left[ \min \left\{ B_1, u_1(\x) \right\} \, / \, 
\min \left\{ B_1, v_{10}(\x) \right\} \right]$ over $\domain$.
Based on the experimental evidence, the resulting efficiency 
gains are quite substantial.  E.g., in the examples of sections 
\ref{ss:convergence} and \ref{sss:flr_isotropic}, 
the MFL allowed additional exclusion of 
50\% and 93\% of the remaining gridpoints respectively.
In the example of section \ref{ss:two_secondary}, where $r=2$,
the corresponding reduction was 76\%, confirming that
the recursive procedure to compute the MFL when $r>1$
is worthwhile.
\label{rem:domain_reduction_efficiency}
\end{remark}

\subsection{Selecting $\tau$ \& optimal ordering}
\label{ss:budget_delta}

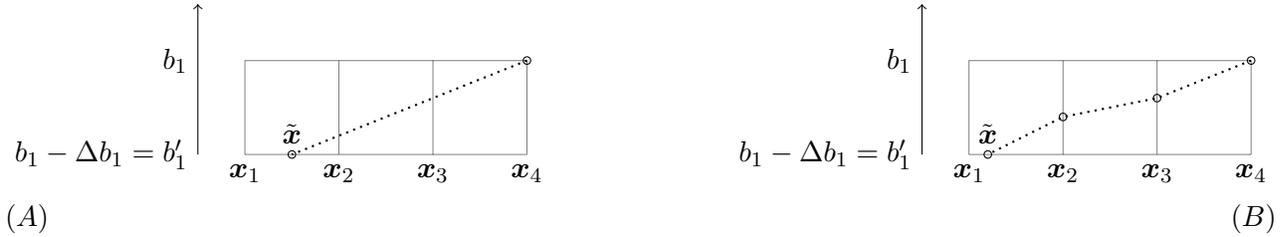
\begin{figure}[hhhh]
$
\hspace*{-2cm}
\begin{array}{lr}
\begin{tikzpicture}[scale=2.5]

	\draw[->]		(0.25,0) node[left]{$b_1 - \Delta b_1 = b_1'$} -- (0.25,0.5)	node[left]{$b_1$} -- (0.25,0.8);

	\draw[xstep=.5cm,ystep=.5cm,gray,very thin] (0.5,0) grid (2,0.5);

	\draw			(0.5,0) node[below]{$\x_1$};
	\draw			(1.0,0) node[below]{$\x_2$};
	\draw			(1.5,0) node[below]{$\x_3$};
	\draw			(2.0,0) node[below]{$\x_4$};

	\draw [thick, dotted] (2.0, 0.5) -- (0.75, 0) node[above]{$\xtilde$};
	\draw			(0.75, 0.0) circle (.2mm);
	\draw			(2.0, 0.5) circle (.2mm);

\end{tikzpicture}&
\hspace*{2cm}
\begin{tikzpicture}[scale=2.5]

	\draw[->]		(0.25,0) node[left]{$b_1 - \Delta b_1 = b_1'$} -- (0.25,0.5)	node[left]{$b_1$} -- (0.25,0.8);

	\draw[xstep=.5cm,ystep=.5cm,gray,very thin] (0.5,0) grid (2,0.5);

	\draw			(0.5,0) node[below]{$\x_1$};
	\draw			(1.0,0) node[below]{$\x_2$};
	\draw			(1.5,0) node[below]{$\x_3$};
	\draw			(2.0,0) node[below]{$\x_4$};

	\draw [thick, dotted] (2.0, 0.5) -- (1.5, 0.3) -- 
						  (1.0, 0.2) -- (0.6, 0)  node[above]{$\xtilde$};
	\draw			(0.6, 0.0) circle (.2mm);
	\draw			(1.0, 0.2) circle (.2mm);
	\draw			(1.5, 0.3) circle (.2mm);
	\draw			(2.0, 0.5) circle (.2mm);
\end{tikzpicture}\\
(A) & (B)
\end{array}
$
\caption{
{\footnotesize
Choices of $\tau$ illustrated for $r=1$ and $n=1$.
Let $\xhat = (\x_4, b_1)$.
A trajectory corresponding to a particular $\ba$ is shown
by a dotted line.  
Algorithm 1 is illustrated on the left: set
$\tau_{\blds{a}}  = \Delta {b_1} / K_1(\x, \ba)$ and assume
that $\fB$ and $K_i$'s don't change.  This ensures
$\tilde{b}_1 = b_1'$ and gives a direct formula for $\xtilde \in [\x_1, \x_2].$
Algorithm 2 is illustrated on the right: $\fB$ and $K_i$'s are 
re-evaluated at each intersection with 1-dimensional cells.
Algorithm 3 is based on the fact that, if $W(\x_3, b_1)$ 
has been already computed, then there is no need to continue 
beyond the first intersection.  In that case, a much smaller 
$\tau_{\blds{a}}$ already guarantees the causality, $\xtilde = \x_3$,
and $\tilde{b}_1 \in [b_1', b_1].$
}
}
\label{fig:setting_tau}
\end{figure}

Formula~\eqref{eq:newstate_after_one_step} for the new state
after using control $\ba$ for $\tau_{\blds{a}}$ seconds
is based on the assumption that, for every fixed $\ba$, 
$\fB$ and all $K_i$'s are approximately constant near $\x$.
There are two advantages to this formula.  First,
it allows to compute $(\xtilde, \btilde)$ very quickly.
Second, taking $\tau_{\blds{a}}  = \Delta {b_1} / K_1(\x, \ba)$
ensures that the new state is in the previous 
$b_1$-slice; i.e., $\tilde{b}_1 = b_1 - \Delta b_1$.
This allows for the explicit causality and marching in the $b_1$
direction.  In section \ref{ss:computational_cost},
we will refer to this implementation as Algorithm 1
(see Figure \ref{fig:setting_tau}A).

There are also two obvious drawbacks to this discretization. 
First, this linear approximation is poor wherever $\fB$ and 
$K_i$'s vary significantly near $\x$.  
Second, the local error in formula~\eqref{eq:discr_augmented} is 
%      COMMENT REMOVED
generally $O(\tau_{\blds{a}}^2)$ even for a smooth $w$ and perfectly known
$W(\xtilde, \btilde)$; thus, it would be preferable
to select the smallest $\tau_{\blds{a}}$ that still allows 
explicit marching in the direction $b_1$.

To address the first of these drawbacks, we have implemented 
Algorithm 2.
Given the current state $\xhat = (\x, \bb)$
and a particular control value $\ba$, we start with the ray from $\xhat$
in the direction $[\fB(\x,\ba), -K_1(\x,\ba), \ldots, -K_r(\x, \ba)].$
We find the first intersection of that ray with
an $(n+r-1)$-dimensional cell of the grid.  If that cell is
fully in the previous $b_1$-slice, we can interpolate, as in Algorithm 1.
Otherwise, we re-evaluate $\fB$, and $K_i$'s at the new point
and follow the new ray, repeating the process until we reach 
the previous $b_1$-slice (see Figure \ref{fig:setting_tau}B).
This procedure is computationally more expensive than Algorithm 1,
since we might need to traverse through several 
$(n+r)$-dimensional cells before finding $\xtilde$,
especially when $\Delta b_1 / min (K_1) \; > \; h / min |\fB|$.
However, it is much more suitable for problems with 
rapidly varying $\fB$ and $K_i$'s.  
(Indeed, in section \ref{ss:computational_cost} we consider several numerical examples,
where these functions are actually discontinuous in $\x$.)

To alleviate the second problem, we have implemented
Algorithm 3.  The idea is to select the smallest
$\tau_{\blds{a}}$ needed for the causality.
Our implementation uses the smallest $\tau_{\blds{a}}$ sufficient 
to guarantee that $(\xtilde, \btilde)$ lies in 
an $(n+r-1)$-dimensional cell, all vertices of which have 
already had values of $W$ previously computed 
(making interpolation possible).
This will certainly be the case if that cell lies fully
in the previous $b_1$-slice (as in Algorithm 2),
but it may also happen earlier;
see Figure \ref{fig:setting_tau}B for an example.
We note that all relevant intersection points are already computed
in Algorithm 2.  If one of them is suitable in the above sense,
we don't need to continue beyond it 
(thus reducing both the computational cost and 
the local truncation error).  In the worst case, we still need 
to trace this piecewise-linear trajectory up to its intersection
with the previous $b_1$-slice (as in Algorithm 2).
As a result, the computational cost of this Algorithm is always
somewhat less than that of Algorithm 2 though obviously higher 
than that of Algorithm 1; see section \ref{ss:computational_cost}.

\begin{remark}
We note that both the accuracy and the efficiency 
of Algorithm 3 clearly depend on the order
of computing $W$'s {\em within} the same $b_1$-slice.  
(The location of $\xtilde$ in Figure \ref{fig:setting_tau}B will be different
depending on whether we have already computed $W(\x_3, \bb_1)$
prior to $W(\x_4, \bb_1)$.)
Our current implementation uses a simple {\em lexicographic ordering}
in secondary costs and then in spatial directions:
assuming that $\xhat = (\x, \bb) = (x_1, \ldots, x_n, b_1, \ldots, b_r)$,
we compute within each $b_1$-slice in the direction $b_2$,
within each $(b_1,b_2)$-slice in the direction $b_3$, $\; ... \;$ , 
within each $(\bb)$-slice in the direction $x_1$, 
within each $(\bb, x_1)$-slice in the direction $x_2$, etc.
While the preference for the secondary-cost (over spatial) directions is clearly 
motivated by the explicit causality (i.e., $K_i$'s are positive), 
our fixed ordering of the spatial directions is arbitrary and hardly optimal.
In the future, we would like to investigate the effect of using 
different orderings in spatial directions
(e.g., in the $r=1$ case, sorting $\xhat$'s
by $W$ values found in the previous $b_1$-slice).

We emphasize that the goal of finding a good ordering within a 
$b_1$-slice is simply to reduce the local truncation error and
to speed up each $b_1$-step.  
\label{rem:alg3_ordering}
\end{remark}

\section{Numerical Experiments.}
\label{s:numerical-results}

We illustrate our numerical method with several examples. 
Our approach requires to choose a 
primary objective function and treat other objectives as secondary.
All feasible controls will satisfy the constraints $\J_i \leq B_i$ for $i=1, \ldots, r.$ 
We then minimize $\J_{0}$ along these feasible paths.

All of the examples in this section are computed for a two-dimensional system state 
and we assume that $\cdomain = [0,1] \times [0,1].$
In all of our examples (except the convergence study in section \ref{ss:convergence}), 
the goal is to reach a single target point.
We will thus assume that
$q_i=0$ at the target and $q_i=+\infty$ on the rest of $\boundary$ for $i=0, \ldots, r.$
Wherever not explicitly specified, the target is at the point $(1,1)$.

Most of the examples (except for section \ref{sss:flr_anisotropic}) are fully
isotropic in cost \& dynamics, and the results are obtained
by solving a variant of equation~\eqref{eq:multiobject_Eikonal}. 
For the isotropic case, the dynamics of the system simplifies as
$\y'(t) = f(\y(t)) \begin{bmatrix} \cos\left(\alpha(t)\right) \\\sin(\alpha(t))\end{bmatrix},$
where $\y(0) = \x \in \domain \subset \R^2$ and $f: \cdomain \mapsto \R$ is the speed,
and $\alpha \in \left[0,2\pi\right]$ is the control parameter.

In each test problem, we plot constrained-optimal paths 
for different resource vectors $\bb$.  We also show {\em secondary-optimal trajectories}
(recovered from $u_i$'s, for $i=1,  \ldots, r$) even in cases when these trajectories do not 
satisfy other integral constraints.

We note that the observability examples in subsections \ref{ss:hiding} - \ref{ss:two_secondary}
use piecewise continuous running costs and/or speed functions.
Even though this violates assumption (A1), the value function is also
well-defined in this case and the semi-Lagrangian discretization
(based on discretizing Bellman's optimality principle) appear to
converge to it correctly.  
A detailed discussion of viscosity solutions to HJB PDEs with
discontinuous Lagrangian can be found in \cite{Soravia_DiscontinLagrangian}.
A different class of observability-influenced path planning problems 
is considered in \cite{ChengTsai}.

All Figures presented in this section were computed
using Algorithm 1 (in subsections 
\ref{ss:convergence} - \ref{ss:flight-path},
where the running costs and speeds are continuous)
and Algorithm 3 (in subsections \ref{ss:hiding} - \ref{ss:two_secondary}).

\subsection{A simple example: convergence study.}
\label{ss:convergence}

\begin{figure}[hhhh]
\centerline{
$
\begin{array}{cc}
\begin{array}{c}
\psfig{file=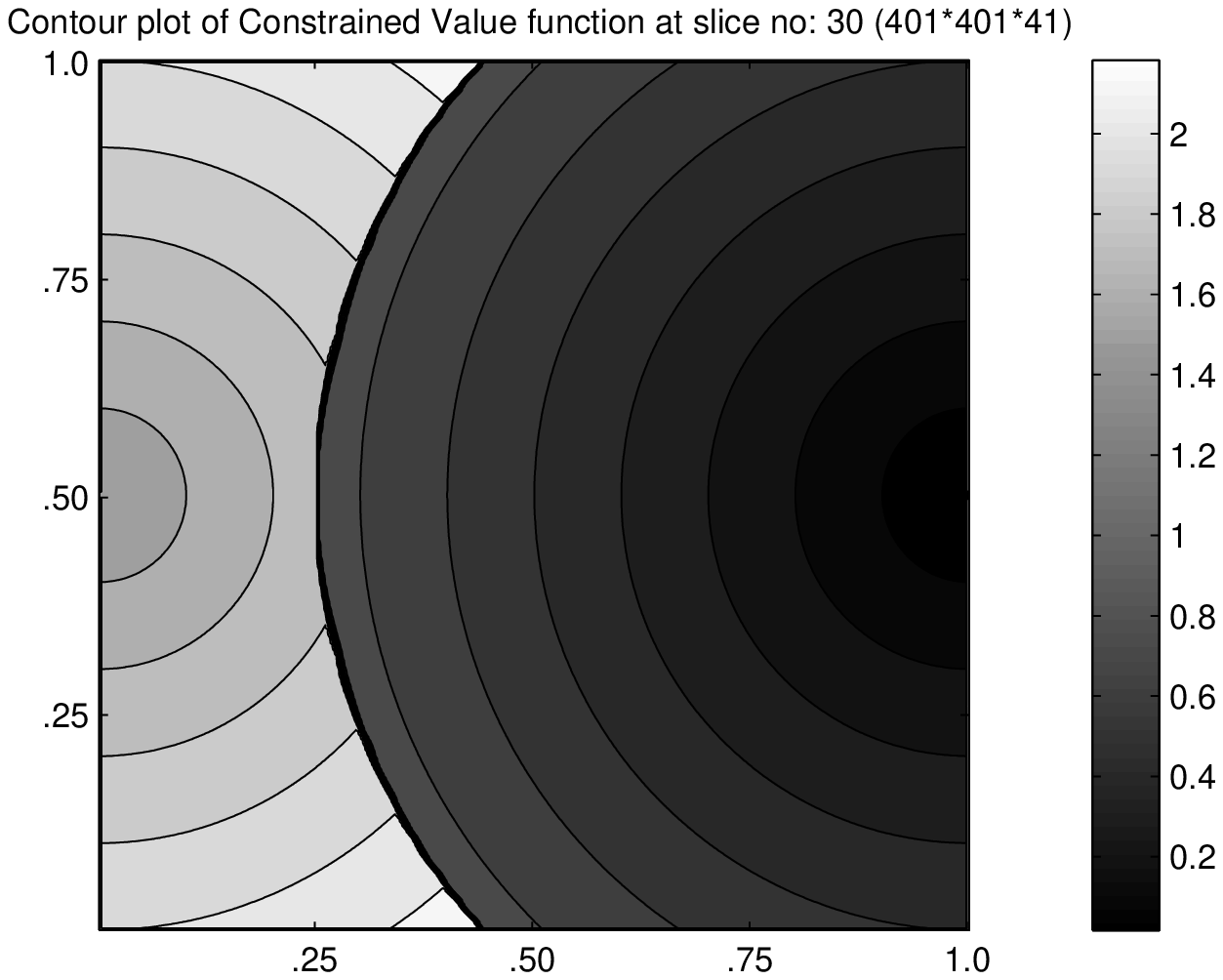,height=8cm}\\
\hspace{.2in}\\
\end{array}&
\begin{array}{c}
     \begin{tabular}{|l|c|c|c|} \hline
	 \multicolumn{4}{|c|}{$L_1$ error computed on $\domain_f$}\\
	 \hline
     \backslashbox{$h$}{$\Delta b_1$}&
     $1/10$ &
     $1/20$ &
     $1/40$\\
    \hline
    {\bf $1/40$ }		& 0.0161431 & 0.0132897 & 0.0217214\\
    \hline
    {\bf $1/80$ }		& 0.0117587 & 0.0078472 & 0.0081641\\
    \hline
    {\bf $1/160$ }		& 0.0088235 & 0.0051967 & 0.0044650\\
    \hline
    {\bf $1/320$ }		& 0.0075405 & 0.0034618 & 0.0024226\\
    \hline
    {\bf $1/640$ }		& 0.0069347 & 0.0026052 & 0.0014743\\
    \hline
     \end{tabular}
\\
\\
     \begin{tabular}{|l|c|c|c|} \hline
	 \multicolumn{4}{|c|}{$L_{\infty}$ error computed on $\domain'(\Delta b_1)$}\\
	 \hline
     \backslashbox{$h$}{$\Delta b_1$}&
     $1/10$ &
     $1/20$ &
     $1/40$\\
    \hline
    {\bf $1/40$ }		& 0.0039922 & 0.0114681 & 0.3558851\\
    \hline
    {\bf $1/80$ }		& 0.0016614 & 0.0024012 & 0.0533905\\
    \hline
    {\bf $1/160$ }		& 0.0007978 & 0.0008307 & 0.0038505\\
    \hline
    {\bf $1/320$ }		& 0.0004434 & 0.0003989 & 0.0004491\\
    \hline
    {\bf $1/640$ }		& 0.0002202 & 0.0002217 & 0.0001994\\
    \hline
     \end{tabular}
\\
\end{array}\\
%      COMMENT REMOVED
\end{array}
$
}
\caption{
\footnotesize{
Convergence to a discontinuous solution (section \ref{ss:convergence}). Left: the level sets of $w(x, 0.75)$.
Right: the $L_1$  and $L_{\infty}$ errors computed for various $(h, \Delta b_1)$ combinations. 
}}
\label{fig:convergence}
\end{figure}

We use a two-dimensional generalization of example in Figure \ref{fig:secondary_ex1}
to test the convergence of our method.
We assume that $f = 1$, and $K_0 = K_1 = 1$.
Unlike in all other examples of this section, here we assume 
that there are two possible exit points on the boundary
$A_1 = (0,0.5)$ and $A_2 = (1, 0.5)$ and define
$$
\begin{array}{ccc}
q_0(\x) = 
\begin{cases}
1.5, & \text{ if } \x = A_1;\\
0, & \text{ if } \x = A_2;\\
\infty, & \text{ otherwise.}
\end{cases} 
& \qquad &
q_1(\x) = 
\begin{cases}
0, & \text{ if } \x = A_1;\\
0, & \text{ if } \x = A_2;\\
\infty, & \text{ otherwise.}
\end{cases}
\end{array}
$$ 
As a result, $\J_1$ is simply the pathlength of the chosen trajectory,
all (constrained optimal) trajectories are straight lines leading to $A_1$ or $A_2$,
and the analytical expression for the discontinuous solution $w$ is readily available:
$$
w(\x,b_1) = 
\begin{cases}
| \x - A_2 |,			& \text{ if } | \x - A_2 | \leq b_1;\\
| \x - A_1 | + 1.5,	& \text{ if } | \x - A_1 | \leq b_1 < | \x - A_2 |;\\
\infty, & \text{ otherwise.}
\end{cases}
$$
We use this example to test the convergence of our method numerically.
(See Figure \ref{fig:convergence}.)

We recall that there are typically four different sources of error in 
semi-Lagrangian schemes:  
(1) due to approximating the boundary conditions (in our case --
also the approximation of the MFL) on the grid;
(2) due to interpolation at the foot of the characteristic 
(e.g., as in formula \eqref{eq:bilinear_interp});
(3) due to approximating the characteristic and
(4) due to approximate integration of the solution along that characteristic.

Since in this example the characteristics are straight lines 
and $K_0$ is constant, the last two of these four sources of error are absent.
Of the remaining two, the first should clearly be decreasing with $\Delta b_1$
(since errors of order $O(\Delta b_1)$ are introduced when representing
the MFL on the grid; see Remark \ref{rem:MFL}).  
Assuming that the solution is smooth, each interpolation error is 
$O(h^2)$ (due to our use of bilinear interpolation \eqref{eq:bilinear_interp}).
However, the number of interpolations is inversely proportional to $\Delta b_1$
and the cumulative effect of interpolation errors is larger 
when the total number of $b_1$-slices increases. 
Thus, if $h$ is held constant, decreasing $\Delta b_1$ will eventually 
result in an increase in the overall error;  this can be seen in the top 
two rows of the corresponding $L_1$-errors table.  Of course, in more general problems,
where characteristics are not straight lines and $K_i$'s are not constant,
the third and fourth sources of errors would normally prevent this phenomenon.

To verify the convergence, we include two tables of error measurements
for various $(h, \Delta b_1)$ combinations.
The first table shows errors measured in $L_1$ norm on $\domain_f$ 
(the subset of $\domain_e$ on which $w$ is finite though not necessarily continuous).  
For discontinuous solutions, convergence in $L_{\infty}$ norm 
is generally possible only for numerical methods that 
explicitly track the location of discontinuities.  
Since no such explicit tracking is performed here, we can only demonstrate 
$L_{\infty}$-convergence away from discontinuities. 
The theoretical results in \cite{BardiFalconeSoravia1,BardiBottacinFalcone,BardiFalconeSoravia2} 
guarantee that, if $h = o(\tau)$ as $\tau=\Delta b_1 \rightarrow 0$, then
the semi-Lagrangian schemes uniformly converge to the viscosity solution $w$
on any compact subset on which $w$ is continuous. 
For each value of $b_1$, we define the discontinuity set of $w$ in the 
corresponding $b_1$-slice:  
$$
\mathcal{D}_{b_1} = \left\{ (\x, b_1) \in \domain_e \, \mid \, 
w \text{ is discontinuous at } (\x,b_1) \right\},
$$
and for each $\Delta b_1$ we define a subset of $\domain_e$
to study the convergence:
$$
\domain'(\Delta b_1) = \left\{ 
\xhat = (\x, b_1) \in \domain_f \mid
\distance ( \x, \mathcal{D}_{b_1} ) \, \geq \, \frac{3}{2} \Delta b_1 
\right\}.
$$
The second table shows $L_{\infty}$ errors measured on $\domain'(\Delta b_1)$.
Since the width of the ``excluded band'' is defined for each column of the second table separately,
the initial errors in the second and third columns are significantly larger, 
but quickly decrease as $h \rightarrow 0$.
In this example, our numerical results suggest the convergence 
even for $h = O(\tau)$.

\subsection{Fastest paths (with restriction on pathlength).}
\label{ss:fastest_length_restricted}

Here we consider two different examples (one isotropic with obstacles, 
the other anisotropic without obstacles, both inhomogeneous), in which the goal
is to minimize the time-to-target subject to constraints on the maximal
allowable pathlength.  
%      COMMENT REMOVED
In the absence of obstacles, we use $K_0(\x)=1$ and 
$K_1(\x, \ba) = |\fB(\x, \ba)|$ to ensure that 
$\J_0$ and $\J_1$ are respectively the time and the pathlength 
along the corresponding trajectory.  If obstacles are present,
we set $K_0 = +\infty$ inside them, to ensure that all trajectories
passing through them have infinite cost $\J_0$.

\subsubsection{Isotropic dynamics/cost in the presence of obstacles.} 
\label{sss:flr_isotropic}
Here we suppose that the dynamics is isotropic; i.e.,
$\fB ( \x, \ba ) = f( \x ) \ba$, where $\ba \in S_1 = A$.
For every $\x = (x,y) \in \domain$ outside of obstacles, we
will assume that
$f(x,y) \; = \; 1 \, + \, 0.8 \, \sin(4 \pi x) \, \sin(6 \pi y)$
and $K_1(x,y) = f(x,y)$.
\begin{figure}[hhhh]
%      COMMENT REMOVED
\epsfig{file=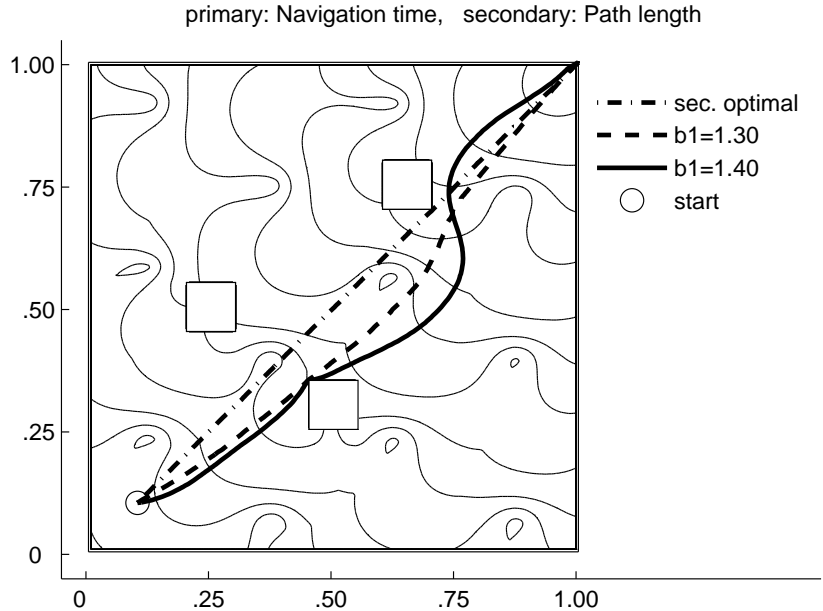,width=13cm}
\caption{
Fastest paths, isotropic dynamics (constrained by pathlength) and obstacles.
Optimal trajectories and the level curves of $u_0$.}
\label{fig:fastest_length_restricted}
\end{figure}
Consider a goal of finding the fastest path to the top right corner
$(1,1)$ subject to restriction on a maximum allowable pathlength
and in the presence of three rectangular obstacles.
The computations are performed on a $301 \times 301 \times 301$
grid and several trajectories are shown in Figure \ref{fig:fastest_length_restricted}
superimposed on top of the level curves of $u_0$.
Secondary-optimal trajectory (dotted line) is the shortest path.
Bold solid line shows a trajectory computed for a large
maximum pathlength.  The constraint is slack in this case and 
the resulting trajectory is in fact primary-optimal; 
hence its orthogonality to the level curves of $u_0$ 
(since for unconstrained isotropic problems the characteristics coincide
with the gradient lines of $u_0$).
The dashed line shows a constrained-optimal trajectory 
for a smaller budget (binding constraint).

\subsubsection{Anisotropic inhomogeneous dynamics.} 
\label{sss:flr_anisotropic}

The following example of anisotropic dynamics was previously used in
\cite{SethVlad2, SethVlad3}, where it was motivated by problems in anisotropic seismic imaging.
This belongs to a class of anisotropic geometric dynamics problems, i.e.,
$\fB ( \x, \ba ) = f( \x, \ba ) \ba$, where $\ba \in S_1 = A$.
Suppose $C:[0,1] \mapsto \R$ is a smooth function.  We are interested in defining $f$ 
so that, for every $\x = (x,y) \in \domain$, the ``velocity profile'' 
$V(\x) = \{f( \x, \ba ) \ba \, \mid \, \ba \in S_1 \}$ is an ellipse
whose major/minor semi-axis have lengths $F_2$ and $F_1$ respectively and the major
semi-axis is aligned with the graph of $C$ (i.e., parallel to its tangent at the point $x$).
This is attained by setting
$$
f( \x, \ba ) = F_2 
\left(
1 +  \left( \left[
\begin{array}{c}
p \\ q \\
\end{array}
\right] \cdot \ba \right)^2
\right)^{-1/2},
\qquad
\text{ where }
\left [
\begin{array}{c}
p \\ q\\
\end{array}
\right ] = \frac { \sqrt{\left(\frac{F_2}{F_1}\right)^2 -1 } }
{ \sqrt{ 1 + \left(\frac{dC}{dx}(x)\right)^2 }}
\left [
\begin{array}{c}
\frac{dC}{dx}(x) \\ -1 \\
\end{array}
\right ].
$$
\begin{figure}[hhhh]
$
\hspace*{-14mm}
\begin{array}{cc}
\epsfig{file=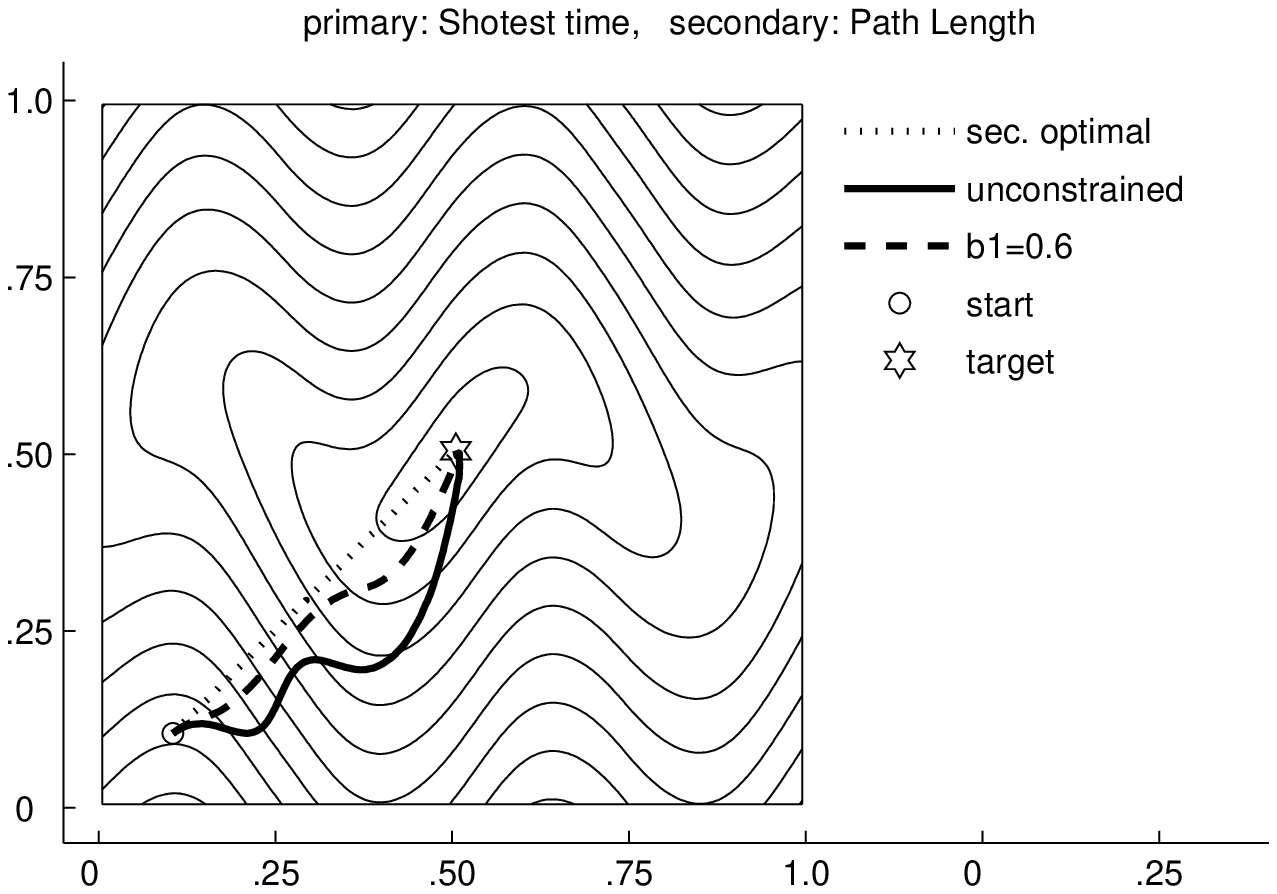,height=7cm}&
\epsfig{file=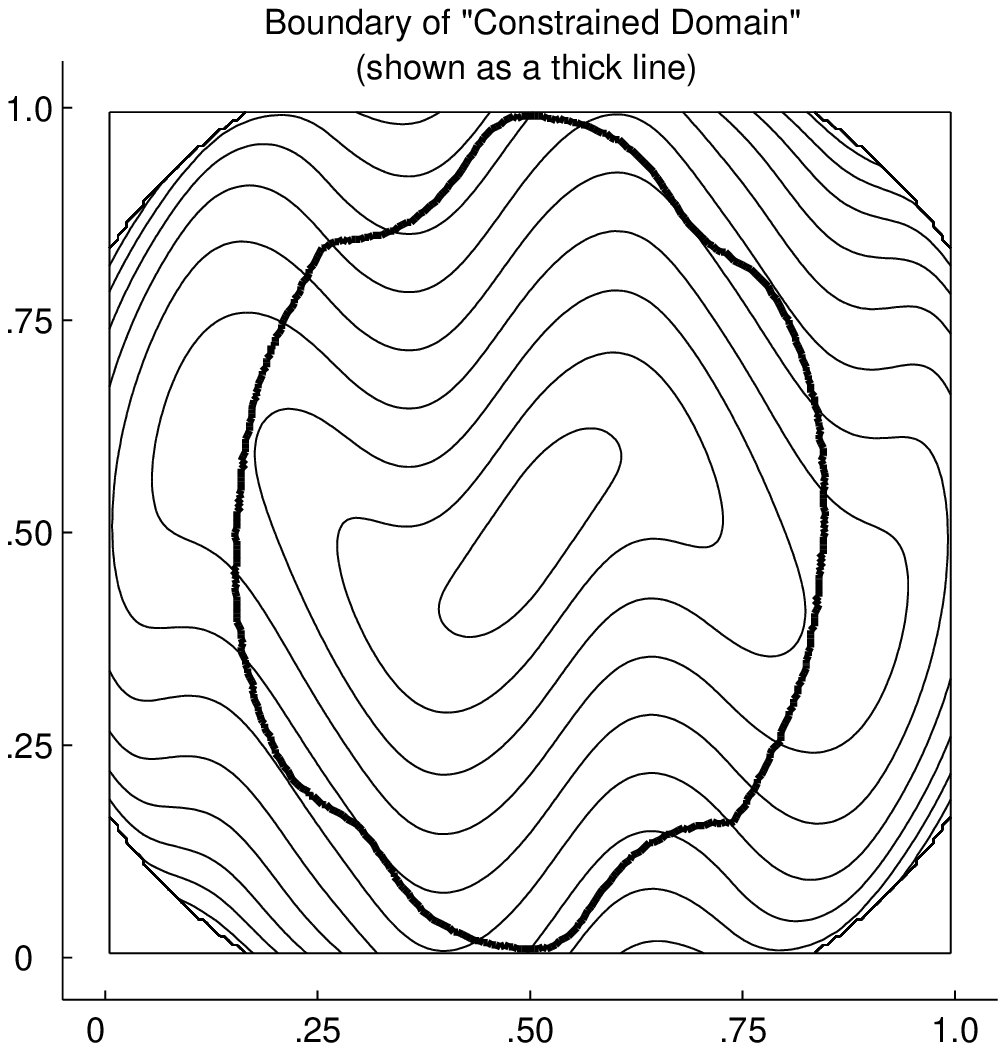,height=6.5cm}\\
{\bf A} & {\bf B}
\end{array}
$
\caption{
Anisotropic dynamics, fastest paths (constrained by pathlength).  
Left: level sets of $w(\x, 1.5)$ (corresponding to level sets of $u_0$ in 
Figure 6A of \cite{SethVlad2}) and optimal trajectories from $(0.1, 0.1)$
to $(0.5, 0.5)$.  Right:  level sets of $w(\x, 0.6)$.
}
\label{fig:anisotropic_example}
\end{figure}
The derivation can be found in \cite{SethVlad3} and \cite{VladThesis}. 
We use $C(x) = 0.1225 \sin(4\pi x)$, $F_2 = 0.8$, $F_1 =0.2$ and compute 
the optimal (fastest) trajectories to the center of $\domain$.
The unconstrained problem leads to an anisotropic static Hamilton-Jacobi PDE
for $u_0$, which can be efficiently solved using Ordered Upwind Methods; 
the level sets of $u_0$ for the above parameters can be found in 
Figure 6A of \cite{SethVlad2}).  
Here we compute pathlength-constrained min-time trajectories for the same example
on a $201 \times 201 \times 301$ grid.  Figure \ref{fig:anisotropic_example}A
shows optimal trajectories to the center starting from the point $(0.1, 0.1)$
superimposed on the level sets of $w(\x, 1.5)$ .  This is a large ``budget'': 
since $b_1 = 1.5 \, > \, v_{10}(\x)$ for $\forall \x \in \domain$, 
we have $w(\x, 1.5) = u_0(\x)$ and all unconstrained-optimal trajectories 
are feasible.   
(We note that these optimal trajectories do not follow the 
gradient lines of $u_0$ since the speed $f( \x, \ba )$ is anisotropic.)
On the other hand, for $b_1 = 0.6$, parts of the domain
are not reachable (note that $w = +\infty$ in the corners of the domain in 
Figure \ref{fig:anisotropic_example}B).  
The thick line separates the part of $\domain$ where $w(\x, 0.6) = u_0(\x)$;
since the starting point $(0.1, 0.1)$ is outside that set, 
its constrained-optimal trajectory is different from the unconstrained-optimal.

\subsection{Optimal flight-path: minimizing weather threat.}
\label{ss:flight-path}
Here, we consider a test example introduced by Mitchell and Sastry \cite{MitchellSastry}: 
finding an optimal path for an airplane flying from one city to another. 
In  one of their examples, weather threat is the primary objective to be minimized	
while the fuel minimization is the secondary objective. The running fuel cost $K_{1}$ and the speed of 
the airplane $f$ are both taken to be unity in this example. 
Hence, the secondary objective, fuel cost is same as the path length.

\begin{figure}[hhhh]
\epsfig{file=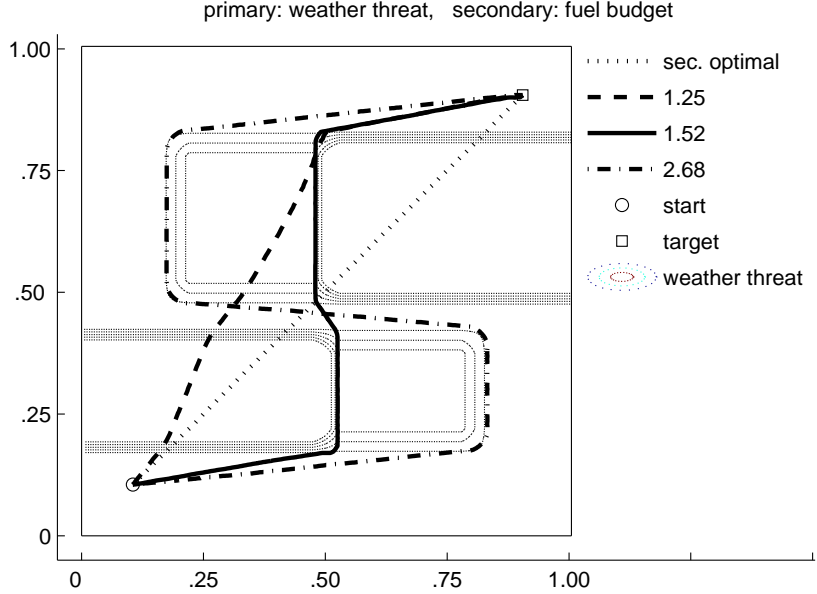,width=13cm}
\caption{Path planning for an airplane striving to be away from the region 
having high weather threat probability.
The optimal paths are shown for the same constraint values used
in Figure 3 of \cite{MitchellSastry}.}
\label{fig:mitchell}
\end{figure}

As mentioned in \cite{MitchellSastry}, weather threat at a location, intuitively, is 
the probability of encountering a storm.  Assuming the probability at all locations to be independent, 
the total probability of encountering a storm during the flight is one minus the product of 
the probabilities of not encountering the storm at all locations along the path. 
As we require running costs of integral type, 
the logarithm of these probabilities are used to make weather-threat cost additive in nature. 
The numerical parameters for this simulation are listed in table 1. 
Figure~\ref{fig:mitchell} shows optimal paths corresponding to different bounds on fuel budget, 
as well as the contours of weather threat cost. The latter is taken to be unity everywhere in $\cdomain$ 
except within the two rectangular bars where it jumps to a relatively higher value discontinuously. 
The magnitude of weather threat is 12 in the darker part of each of the rectangular bars and 
4 in the other part of the bar. Weather threat function is further smoothened to remove 
the discontinuity and make it Lipschitz continuous.

As expected, the secondary-optimal path (the shortest path) follows a straight line between the two cities. 
Bounds on fuel budget associated with other constrained-optimal paths are shown in the legend. 
The Figure shows that, as the fuel budget increases, the airplane chooses a path through 
the region less susceptible to weather threat. These optimal paths match closely those 
in \cite{MitchellSastry}.

\begin{figure}[hhhh]
\epsfig{file=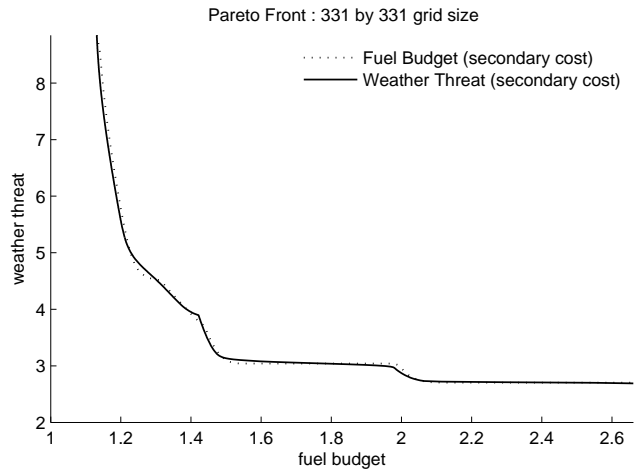,width=10cm}
\caption{ Pareto front for a particular destination for the airplane problem.
A convexification of this front can be found in Figure 4 in \cite{MitchellSastry}.}
\label{fig:pareto}
\end{figure}

Figure~\ref{fig:pareto} shows the Pareto front for this example. 
Pareto front is generated using two different approaches.
In the first approach, fuel budget is taken as secondary objective while 
in the other approach weather threat is the secondary objective. 
As $h$ and $\Delta b_1$ decrease, these two versions of Pareto front
look more and more similar.
%      COMMENT REMOVED
The non-convex parts of this front have not been found in \cite{MitchellSastry};
this illustrates the advantage of our approach as compared to the ``weighted-sum scalarization''
based techniques.

{\footnotesize
\begin{table}[hhhh]
\centering
\begin{tabular}{l l}
 \hline
Grid points in system-domain       &: 201 $\times$ 201\\
Grid points along secondary budget &: 401\\
Speed of the vehicle &: 1\\
Primary cost  (Weather threat)&: explained in 4.1\\
Secondary cost(Fuel rate)     &: 1\\
Starting point for optimal path &: (0.1, 0.1)\\
Target point for optimal path &: (0.9, 0.9)\\
 \hline
\end{tabular}
\caption{Numerical parameters for the optimal flight-path example (Figure \ref{fig:mitchell}).}
\end{table}
}

\subsection{Computation of non-visible region}
\label{ss:non-visible region}
The remaining examples deal with robotic navigation in domains with obstacles
in the presence of friendly and adversarial observers.
The observer's position is assumed to be known and static
and $\domain$ is split into parts directly visible ($\domain_v$) or 
invisible ($\domain_i$) to each observer due to obstacles.   
The observability running cost is then set to be high on 
$\domain_v$ and low on $\domain_i$ for an enemy observer 
(and the opposite of this for a friendly observer).

To enable the computation of constrained-optimal paths,
we need to first determine $\domain_v$ and $\domain_i$.
While this can be accomplished by many methods, we
use the efficient technique based on a 
Fast Marching Method and described in \cite{SethBook2}.

We first solve the Eikonal equation 
$|\grad \psi_1(\x)| f(\x) = 1$ with $f=1$ on $\cdomain$ 
and the boundary condition
$\psi_1=0$ at the observer's location. 
We then solve the Eikonal PDE for $\psi_{2}$ with the same boundary condition
but with $f=0$ inside the obstacles. 
The region with $\psi_{2} > \psi_{1}$ defines the non-visible region. 
In practice, we use a (heuristically adjusted) threshold on the difference of 
$\psi_2$ and $\psi_1$. 
Thick grey lines are used in the following Figures to show the boundaries of $\domain_i$.

\subsection{Avoiding the observer.}
\label{ss:hiding}
Here, we find optimal paths for a robot navigating in a region containing stationary obstacles 
and a stationary enemy observer. The robot strives to be minimally exposed to the enemy at the same 
time making sure to avoid the obstacles and stay within the specified fuel budget.
Observability by the enemy becomes the primary cost ($\J_0$) to be minimized. 
As long as the secondary running cost $K_1$ remains positive, our numerical method can
solve the PDE~\eqref{eq:multiobject_PDE} efficiently even with
with $K_0=0$ observability cost in $\domain_i$.  However, this leads to 
infinitely many possible optimal paths
since any portion of the path inside $\domain_i$ would not contribute anything to $\J_0$. 
To remove this arbitrariness and non-uniqueness, we assigned a small non-zero observability cost $K_0=0.1$ 
on $\domain_i$. Table 2 shows the other numerical parameters used in this experiment. 
\begin{figure}[hhhh]
%      COMMENT REMOVED
\epsfig{file=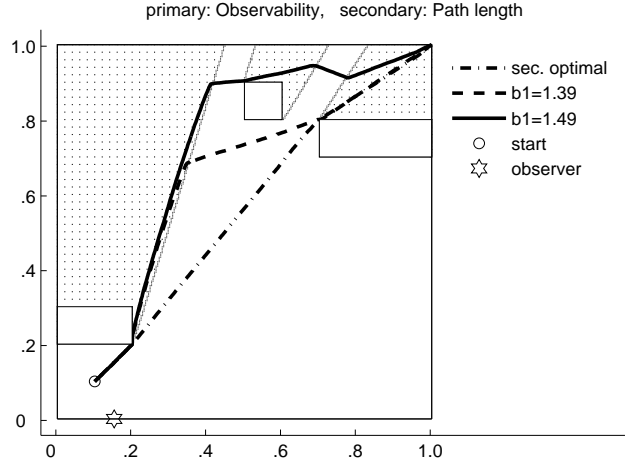,width=10cm}
\caption{Optimal path for a robot navigating to minimize the exposure to a static enemy observer.}
\label{fig:avoiding observer}
\end{figure}
{\footnotesize
\begin{table}[hhhh]
\centering
\begin{tabular}{l l}
 \hline
Grid points in system-domain       &: 251 $\times$ 251\\
Grid points along secondary budget &: 301\\
Speed of the vehicle &: 1\\
Primary cost  (Observability cost)&: 10 \ in the visible region\\
&: 0.1 in the non-visible region\\
Secondary cost(Fuel rate)     &: 1\\
Starting point for optimal path &: (0.1, 0.1)\\
Target point for optimal path &: (1.0, 1.0)\\
Observer location &: (0.15, 0)\\
 \hline
\end{tabular}
\caption{Numerical parameters for the exposure-minimization example (Figure \ref{fig:avoiding observer}).}
\end{table}
}

Figure~\ref{fig:avoiding observer} shows the non-visible region with its boundary as thick solid lines. 
The small rectangles represent obstacles. 
Optimal paths corresponding to different bounds $B_{1}$ on fuel budget are also plotted. 
Since the speed of motion is constant, the running costs are piecewise constant, 
and the obstacles are polygonal, it is easy to prove that all constraint-optimal paths
are piecewise linear.
As expected, the primary-optimal paths creep along the boundary of non-visible region. 
In this example, the target lies in $\domain_i$; 
the robot therefore follows the straight line path 
after it enters that component of $\domain_i$.

\subsection{Minimizing fuel consumption/path length, constrained by enemy observability.}
Here, we find the optimal path for a robot with the same 
two objectives as in the previous section. 
But now the path length (or fuel consumption) is the primary cost 
and the enemy observability is secondary.
The numerical parameters are shown in Table (3). 
\begin{figure}[hhhh]
%      COMMENT REMOVED
\epsfig{file=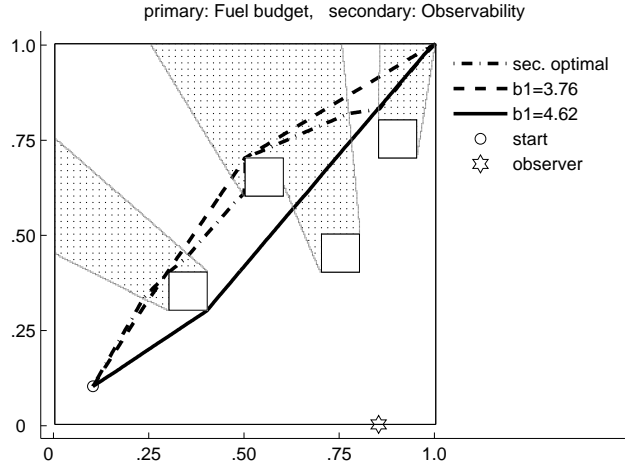,width=10cm}
\caption{Optimal path for a robot minimizing its fuel budget/path length 
under the maximum enemy exposure constraint.  
The enemy exposure budget of 4.62 is large enough to allow 
following the primary-optimal trajectory (thick line).}
\label{fig:shortest path}
\end{figure}
{\footnotesize
\begin{table}[hhhh]
\centering
\begin{tabular}{l l}
 \hline
Grid points in system-domain       &: 301 $\times$ 301\\
Grid points along secondary budget &: 501\\
Speed of the vehicle &: 1\\
Primary cost  (Fuel rate/Path length)&: 1 \\
Secondary cost(Observability)     &: 1 in the non-visible region\\
&: 5 in the visible region\\
Starting point for optimal path &: (0.1, 0.1)\\
Target point for optimal path &: (1.0, 1.0)\\
Observer location &: (0.85, 0)\\
 \hline
\end{tabular}
\caption{Numerical parameters for path length minimization 
constrained by exposure (Figure \ref{fig:shortest path}).}
\end{table}
}
Figure~\ref{fig:shortest path} shows the optimal paths for this example. 
The secondary-optimal path shown as a dashed line is the least exposed to the enemy. 
As expected, the primary-optimal paths become shorter as we 
relax the constraint (or, alternatively, increase the ``budget'') 
for the maximum observability.

\subsection{Striving to be observed}
Given a stationary friendly observer, the goal in many applications is to minimize 
the total time outside of direct visibility while moving to the target.
This is more or less the opposite of the problem considered in section \ref{ss:hiding}.
The total fuel available (or the maximum path length) is still treated as a
constraint. The numerical parameters are shown in table (4).\\

\begin{figure}[hhhh]
%      COMMENT REMOVED
\epsfig{file=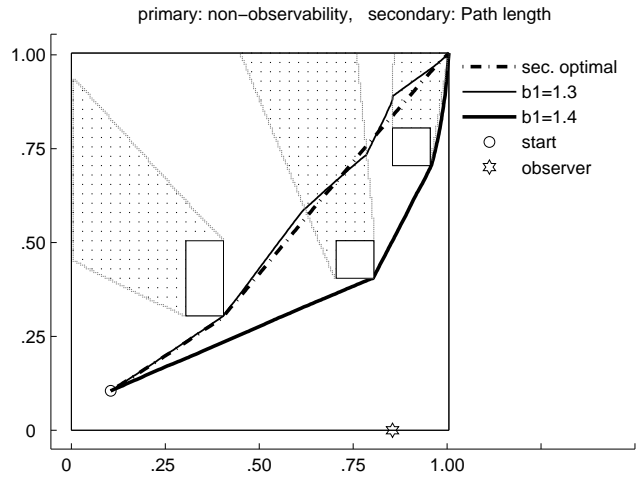,width=10cm}
\caption{ Optimal path for a robot in the presence of a friendly observer.}
\label{fig:non-observability}
\end{figure}

Figure~\ref{fig:non-observability} plots the optimal paths corresponding to different bounds on path length.
When the bound on path length is tight, the robot has no option but to navigate through $\domain_i$.
As we relax this bound, the robot finds a path which is always exposed to the observer.

{\footnotesize
\begin{table}[hhhh]
\centering
\begin{tabular}{l l}
 \hline
Grid points in system-domain       &: 201 $\times$ 201\\
Grid points along secondary budget &: 301\\
Speed of the vehicle &: 1\\
Primary cost(Non-observability)    &: 5 in non-visible region\\
&: 1 in visible region\\
Secondary cost  (Fuel rate)&: 1 \\
Starting point for optimal path &: (0.1, 0.1)\\
Target point for optimal path &: (1.0, 1.0)\\
Observer location &: (0.85, 0)\\
 \hline
\end{tabular}
\caption{Numerical parameters for minimizing the non-exposure
constrained by path length (Figure \ref{fig:non-observability}).}
\end{table}
}

\subsection{Path length minimization subject to two integral constraints.}
\label{ss:two_secondary}
\begin{figure}[hhhh]
$
\begin{array}{c}
%      COMMENT REMOVED
\epsfig{file=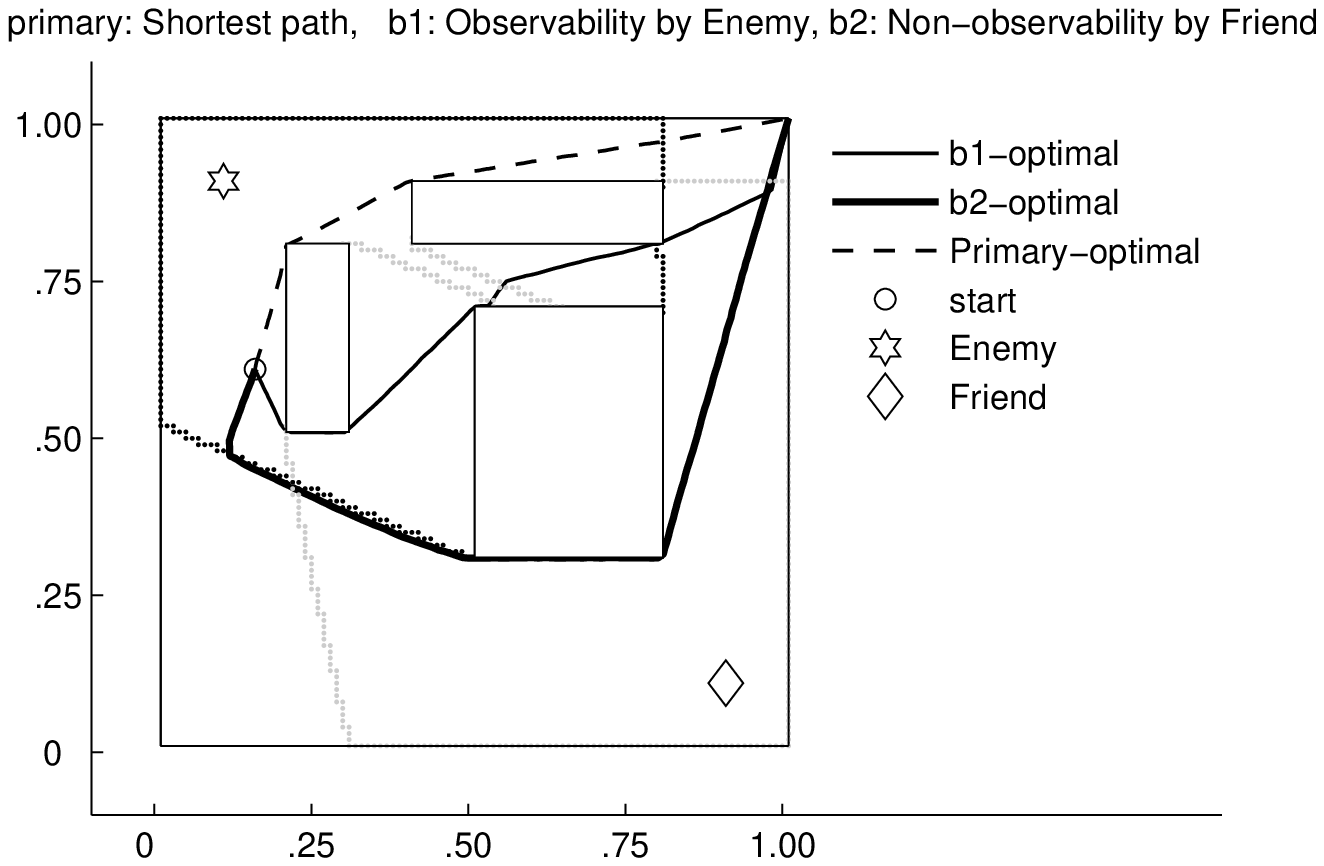,width=14cm}\\
\vspace*{-15mm}
%      COMMENT REMOVED
\epsfig{file=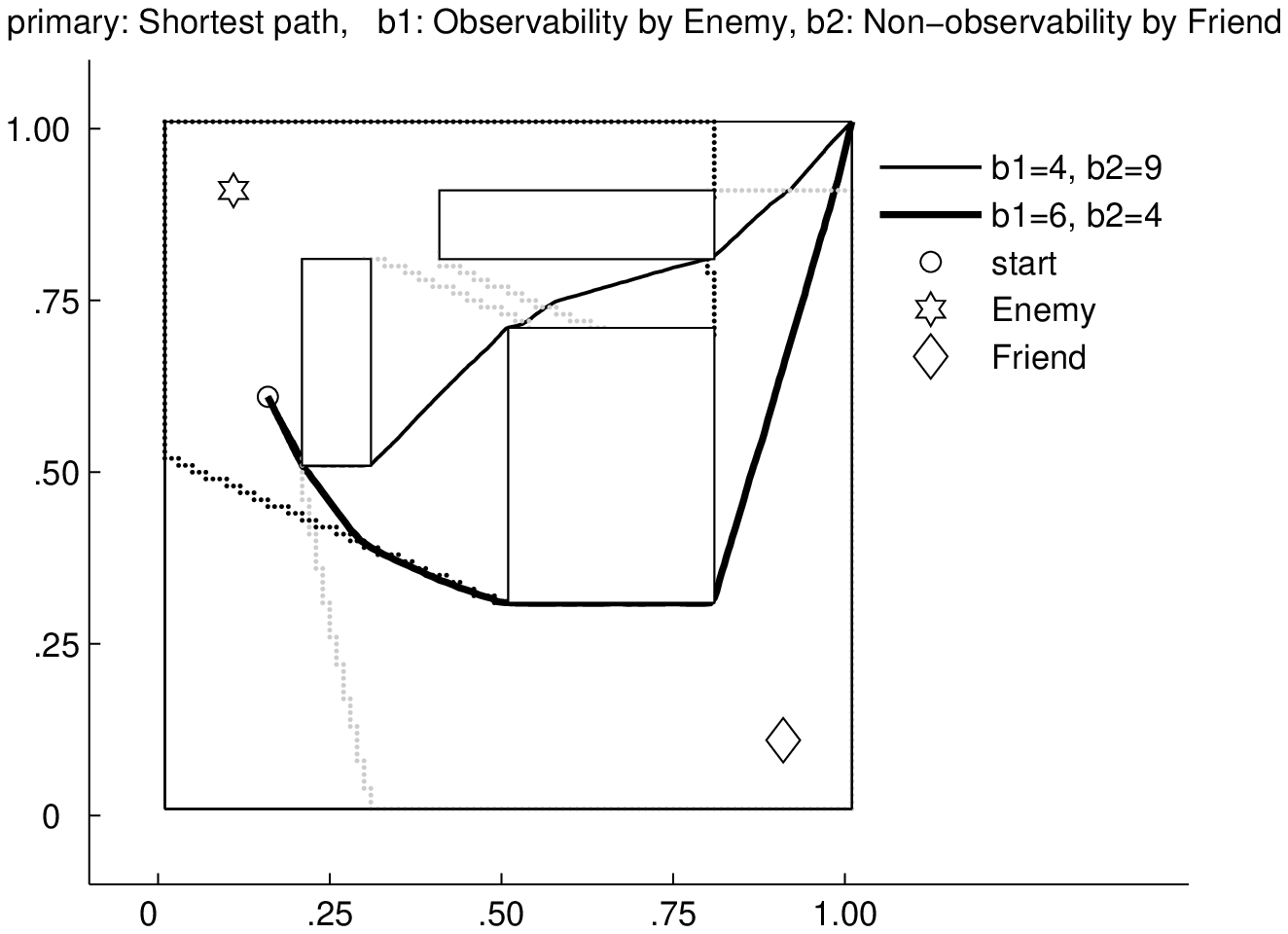,width=14cm}
\end{array}
$
\caption{ A two-secondary costs example: 
fuel-optimality under constraints on visibility by enemy and 
non-visibility by friend.  Dotted lines show the boundaries of 
visibility for both observers. 
Primary and secondary-optimal trajectories (top) and
constrained-optimal trajectories (bottom).
}
\label{fig:two_secondary_costs}
\end{figure}

In this last example, we consider a problem of finding constrained-optimal 
paths in the presence of obstacles and two observers.  The goal is
to minimize the path-length subject to constraints on
the amounts of time the robot can be visible to the enemy observer
and invisible to the friendly observer. Given two secondary costs,
the numerical domain is four-dimensional.  As explained in section
\ref{ss:reduced_domain}, we first solve the PDE ~\eqref{eq:2_cost_init} 
on $\cdomain \times [0,B_2]$
to find the feasibility surface.  We then march in the direction $b_1$
to solve for the value function $w$.  Two constrained-optimal trajectories 
are shown in Figure \ref{fig:two_secondary_costs}.

{\footnotesize
\begin{table}[hhhh]
\centering
\begin{tabular}{l l}
 \hline
Grid points in system-domain       &: 101 $\times$ 101\\
Grid points along each secondary budget &: 301\\
Speed of the vehicle &: 1\\
Primary cost    &: path length\\
Secondary cost 1 &: visibility to enemy \\
Secondary cost 2 &: invisibility to friend \\
Visibility cost values: 1 and 10 for both enemy and observer.\\
\hline
\end{tabular}
\caption{
Numerical parameters for the two-secondary-costs example (Figure \ref{fig:two_secondary_costs}).}
\end{table}
}

\subsection{Discussion of computational complexity.}
\label{ss:computational_cost}
Since our semi-Lagrangian discretization of the PDE~\eqref{eq:multiobject_PDE} 
is explicitly causal, the computational complexity of
the methods is $O(M)$, where $M$ is the number of gridpoint
on $\cdomain_e$.

A careful restriction of the problem to a feasible subdomain
yields an efficient numerical method.  
For the case $r=1$, computations on a 3-dimensional grid
are quite fast on an average laptop.
We have used Dell Inspiron 1505 laptop with 2 GHz Intel Centrino processor, 
and 1 GB RAM.  On a grid with $M=201^3$ gridpoints our
instrumented (and unoptimized) code took 21, 43, and 39 seconds
for Algorithms 1, 2, and 3 of section \ref{s:multi-intro} respectively.
For $r=2$ and $M=101^4$, Algorithm 2 runs in less than 30 minutes on the same laptop.

In principle, only two $b_1$-slices of the grid are needed in RAM 
to enable efficient marching.  However, our current implementation
allocates the entire grid.  As a result, the last example 
(involving $r=2$ and $M=101^2 \times 301^2$ gridpoints and 
$B_1= 11, \, B_2 = 9$) was computed on a machine with $64$ GB RAM
though the memory footprint of the program is $\approx 10$ GB. 
The initialization (by solving for $w_1$ on a $101^2 \times 301$ grid) took 17 seconds.
The Algorithms 1, 2, and 3 then took 27, 67, and 55 minutes respectively.

We note that the computational time is heavily dependent on the values of
$B_1$ and $B_2$.  There are two reasons for this.
First, the tighter constraints make a larger part
of $\cdomain_e$ non-feasible, resulting in a big reduction in 
computational cost for all three algorithms.
Second, in Algorithms 2 and 3, $\tau_{\blds{a}}$ is usually dependent on ratios
between $h$ and $\Delta b_i$'s.  This also influences the number of 
$(n+r)$-dimensional cells traversed 
before reaching a suitable interpolation point $(\xtilde, \btilde)$; 
see section \ref{ss:budget_delta}.  E.g., for $r=1$, when $\Delta b_1 << h$ 
the interpolation is performed after traversing a single cell.
%      COMMENT REMOVED
Decreasing $B_i$ while holding constant the number of gridpoints 
in that direction decreases $\Delta b_i$ proportionally.
To illustrate this point, the problem of section \ref{ss:two_secondary}
on the same grid ($M=101^2 \times 301^2$) ,
but with $B_1 = 4$ and $B_2 =  9$ is much less computationally expensive: 
Algorithms 1,2, and 3 now take only 20, 27, and 26 minutes respectively
on the same computer.

In comparing these execution times to those reported in 
\cite{MitchellSastry}, it is important to keep in mind that 
Mitchell and Sastry have used an upwind finite-difference 
discretization of the Eikonal PDE (i.e., isotropic control problems only).
In contrast, our semi-Lagrangian implementation is also
suitable for much more general (anisotropic
and/or non-small-time-controllable) problems, 
including those, where the minimization in~\eqref{eq:multiobject_PDE}
cannot be performed analytically.

We also note that our current implementation is simple but non-optimal 
since the grid is stored as a multidimensional array
and all the domain-reduction procedures described in section \ref{ss:reduced_domain}
are performed only after the memory for the grid has been allocated.
The ``excluded'' gridpoints are marked (to avoid computing $W$), 
but still take some computational time (in enumerations,
input/output operations, etc.)  A more efficient implementation 
could be clearly built to allocate the memory for non-excluded
gridpoints only, but this would require using non-array 
data structures to represent the grid.

\section{Conclusions.}
\label{s:conclusions}

We have introduced a new numerical method for multiobjective optimal control
and single-objective optimal control in the presence of integral constraints.
Our approach is based on expanding the state space to include 
constraint-budgets and then solving an augmented PDE, whose explicit 
causality allows for a non-iterative (marching) numerical method.
We have also shown the connection between the integral-constrained 
single-objective problem and the task of finding all Pareto-optimal controls.
Our method was illustrated with a number of test-problems for two-dimensional 
optimal control with one and two additional integral constraints.
We have used a flight-path bad weather avoidance example introduced in
\cite{MitchellSastry} as well as several examples of optimal robotic navigation
in the presence of friendly and adversarial observers.

It is a commonly accepted practical rule that lower-dimensional computations
are much less expensive than the higher-dimensional ones (simply because
the number of gridpoints grows exponentially with the dimension).  
However, this simple rule of thumb ignores more subtle issues:
How many PDEs need to be solved on each domain?
How many times do we need to solve each PDE?  
Does the lower-dimensional approach
adequately capture the high-dimensional picture?
Which PDE can be solved with a more efficient numerical method?\\
These questions reflect the difference between our approach and the
prior method by Mitchell and Sastry \cite{MitchellSastry}.
Their method is based on solving a system of $(r+2)$ PDEs
(all but one of them linear, the remaining non-linear PDE
with monotone causality, enabling a Dijkstra-like numerical method)
on $\domain \subset \R^n$ but with an $r$-dimensional parameter space,
which requires solving this system repeatedly $O(2^r)$ times.
Our approach leads to a single PDE with explicit causality
(enabling time-like marching) but on a $(n+r)$-dimensional domain.
Even more importantly, our approach allows recovering the entire 
relevant part of Pareto front, including the non-convex parts of it, which are
inaccessible using the weighted sums method employed in \cite{MitchellSastry}.

The explicit causality of the augmented PDE and 
a careful restriction of the problem to a feasible subdomain
yield an efficient numerical method.  
However, the memory requirements of our method are 
more extensive since at least two $b_1$-slices of the grid have 
to be kept in memory at all times for efficient marching.  
This is an $(n+r-1)$-dimensional grid, in contrast with
an $n$-dimensional grid used by the method in \cite{MitchellSastry}.
Another disadvantage of our approach is the fact that the 
local truncation errors are $O(\hat{h})$ rather than $O(h)$.
As a result, the quality of reconstruction of optimal controls and 
trajectories also depends on $\Delta b_i$'s, 
whereas the method in \cite{MitchellSastry}
can provide a good trajectory reconstruction for each $\lambda$ regardless
of coarseness of the mesh imposed on $\Lambda$.

In the future we would like to build semi-Lagrangian and Eulerian
higher-order accurate methods based on our approach.  
We also intend to explore the use of adaptive grids and unstructured meshes 
(since the constrained-optimal controls can be sensitive
to small changes in available budgets).  It will also be 
relevant to experiment with the efficiency/accuracy implications
of the choice of marching direction --
any other $b_i$ can be chosen since all $K_i$'s are assumed to be positive,
but the time $\tau_{\blds{a}}$ in the semi-Lagrangian
discretization is currently based on $K_1$.

Several other natural extensions
should be possible without breaking the explicit causality 
of the augmented PDE.  We intend to extend our method to
\begin{enumerate}
\item
problems, where $K_i$'s don't have to be positive for $i>1$;
\item
problems with running costs (and exit costs) dependent
on the budgets still remaining;
\item
optimal stochastic control subject to integral constraints;
\item
differential games subject to integral constraints.
\end{enumerate}
The theoretical framework for considering secondary costs
of varying sign has been developed in \cite{Soravia_IntegralConstraint}.
We believe that as long as $K_i >0$ holds at least for one $i\geq 1$, 
the computational efficiency
will not be adversely affected.  If all secondary $K_i$'s 
(but not $K_0$) 
are allowed to change sign, this will require a method based on a 
more subtle monotone causality in $w$.

\vspace{.2in}
\noindent
{\bf{Acknowledgments:}}
The authors would like to thank Ian Mitchell for 
stimulating discussions on different approaches to
continuous multiobjective optimal control. 
The authors are grateful to Maurizio Falcone and 
Roberto Ferretti for suggesting relevant references on
convergence of semi-Lagrangian discretizations.
The anonymous reviewers' comments were also very helpful
in revising this manuscript.

\end{document}